\documentclass{amsart}
\usepackage{amssymb}
\usepackage{amsmath}
\usepackage{latexsym}
\usepackage{eepic}
\usepackage{epsfig}
\usepackage{graphicx}
\usepackage{pb-diagram,lamsarrow,pb-lams,amscd}
\usepackage{eufrak}
\usepackage{mathrsfs}
\usepackage[dutch,english]{babel}

\DeclareMathAlphabet{\mathpzc}{OT1}{pzc}{m}{it}
\newtheorem{theorem}{Theorem}[section]
\newtheorem{corollary}[theorem]{Corollary}
\newtheorem{prop}[theorem]{Proposition}
\newtheorem{lemma}[theorem]{Lemma}
\newtheorem{cor}[theorem]{Corollary}
\newtheorem{definition}[theorem]{Definition}

\newenvironment{remark}{\noindent\textbf{Remark.}}{\qed}

\newcommand{\LL}{\ensuremath{\mathbb{L}}}

\newcommand{\N}{\ensuremath{\mathbb{N}}}
\newcommand{\Z}{\ensuremath{\mathbb{Z}}}
\newcommand{\Q}{\ensuremath{\mathbb{Q}}}

\newcommand{\C}{\ensuremath{\mathbb{C}}}

\newcommand{\OO}{\widehat{\mathcal{O}}}

\newcommand{\A}{\ensuremath{\mathbb{A}}}

\newcommand{\G}{\ensuremath{\mathbb{G}}}

\newcommand{\mX}{\ensuremath{\mathfrak{X}}}
\newcommand{\mY}{\ensuremath{\mathfrak{Y}}}

\newcommand{\cF}{\ensuremath{\mathcal{F}}}

\newcommand{\cY}{\ensuremath{\mathscr{Y}}}

\newcommand{\Spec}{\ensuremath{\mathrm{Spec}\,}}

\renewcommand{\C}{\ensuremath{\mathbb{C}}}

\renewcommand{\A}{\ensuremath{\mathbb{A}}}

\renewcommand{\mX}{\ensuremath{X_\infty}}
\renewcommand{\mY}{\ensuremath{Y_\infty}}
\newcommand{\mZ}{\ensuremath{Z_\infty}}

\renewcommand{\cF}{\ensuremath{\mathscr{F}}}

\renewcommand{\cY}{\ensuremath{\mathscr{Y}}}

\begin{document}
\title{Motivic {S}erre invariants, ramification, and the analytic {M}ilnor fiber}
\author[Johannes Nicaise]{Johannes Nicaise$^{\ast}$}
\thanks{$^{\ast}$ During part of the research for this article, the first author was Research Assistant of the Fund for Scientific Research --
 Flanders (Belgium)(F.W.O.) at the Department of Mathematics of the Katholieke Universiteit Leuven (Belgium).}

\address{Universit\'e Lille 1\\
Laboratoire Painlev\'e, CNRS - UMR 8524\\ Cit\'e Scientifique\\59 655 Villeneuve d'Ascq C\'edex \\
France} \email{johannes.nicaise@math.univ-lille1.fr}


\author{Julien Sebag}
\address{Universit\'e Bordeaux 1,
IMB, Laboratoire A2X, 351 cours de la lib\'eration, 33405 Talence
cedex, France} \email{julien.sebag@math.u-bordeaux1.fr}

 \maketitle
\section{Introduction}
Let us recall the classical definition of a $p$-adic zeta
function, as it was given by Igusa \cite{Igusa:intro}. A survey of
the theory of $p$-adic zeta functions can be found in Denef's
Bourbaki report \cite{DenefBour}.

In its simplest form, a $p$-adic zeta function is defined by the
$p$-adic integral
$$Z_p(f,s)=\int_{\Z_p^{n}}|f(x)|^{s}|dx|\,, $$
where $s$ is a complex variable, $f$ is a polynomial over $\Z_p$
in $n$ variables, $|f|$ is its $p$-adic norm, and $|dx|$ denotes
the Haar measure on the compact group $\Z_p^{n}$, normalized to
give $\Z_p^{n}$ measure $1$. A priori, $Z_p(f,s)$ is only defined
when $\Re(s)>0$. However, Igusa proved, using resolution of
singularities, that it has a meromorphic continuation to the
complex plane. Moreover, it is a rational function in $p^{-s}$.
The numerical data of an embedded resolution for $f$ yield a
complete set of candidate poles of $Z_p(f,s)$, but, since this set
depends on the chosen resolution, a lot of these candidate poles
will not be actual poles of $Z_p(f,s)$. The $p$-adic zeta function
$Z_p(f,s)$ can be written as a Mellin transform of the local
singular series of $f$ (defined by integrating the Gelfand-Leray
form $dx/df$ on the regular fibers of $f$); see
\cite[1.4]{DenefBour}.

When $f$ is defined over $\Q$, Igusa's Monodromy Conjecture
predicts an intriguing connection between the eigenvalues of the
monodromy at complex points of $f^{-1}(0)\subset \C^{n}$, and the
poles of $Z_p(f,s)$, for almost all primes $p$, where $f$ is
considered as a complex, respectively $p$-adic polynomial. Since
the $p$-adic zeta function contains information about the number
of solutions of the congruence $f(x)\equiv 0$ modulo powers of
$p$, the monodromy conjecture establishes a fascinating bridge
between arithmetic properties and complex topology, much like the
Weil conjectures. In fact, we will show in this paper that this
analogy is more than just philosophical.

In order to geometrize the situation, Denef and Loeser introduced
the topological zeta function $Z_{top}(f,s)$, which is some kind
of geometric abstraction of the $p$-adic zeta function. It is a
rational function in the complex variable $s$, defined in terms of
an embedded resolution of singularities, and to show that the
definition does not depend on the choice of resolution, Denef and
Loeser established $Z_{top}$ as a formal limit of $p$-adic zeta
functions. After Kontsevich introduced motivic integration
\cite{Ko}, it became clear that both the $p$-adic and the
topological zeta function are avatars of a universal being, the
``na\"\i ve'' motivic zeta function $Z_{mot}(f,s)$ (the more exact
statement says that $Z_{mot}$ specializes to $Z_{top}$, and to
$Z_p$ for almost all primes $p$). This na\"\i ve motivic zeta
function is defined intrinsically by means of a motivic integral,
but it can be expressed explicitly in terms of an embedded
resolution of singularities for the morphism $f$. In
\cite[3.5.3]{DL3}, Denef and Loeser define the motivic nearby
cycles of $f$ as a limit of the motivic zeta function. It is an
object in an appropriate Grothendieck ring of varieties over the
complex hypersurface $X_s=f^{-1}(0)$, and the fiber over each
point $x\in X_s$ has the same Hodge polynomial as the Milnor fiber
of $f$ at $x$. The motivic monodromy conjecture states that each
pole of the na\"\i ve motivic zeta function induces an eigenvalue
of the monodromy at some point of $X_s$. In this setting, the
notion of pole has to be defined with care \cite[\S 4]{RoVe},
since the Grothendieck ring of varieties is not a domain
\cite{Poo}.

Most research on this conjecture made a ``detour" via a resolution
of singularities, studying the geometry of the exceptional locus
in order to eliminate fake candidate poles, and applying A'Campo's
formula for the monodromy zeta function \cite{A'C}. This approach
has yielded proofs in particular cases, and inspired nice results
concerning the geometry of embedded resolutions, but it seems
difficult to apply this technique to the general case. In this
article, we try to establish a more direct link between the
motivic zeta function of $f$, and the Milnor fibration, using
Berkovich' \'etale cohomology for analytic spaces
\cite{Berk-etale},
 and the theory of motivic integration on rigid spaces
 \cite{motrigid}. More precisely, we introduce the analytic Milnor
 fiber, a rigid variety over $\C((t))$ with the ``same''
 cohomology as the topological Milnor fiber, and whose points are
 closely related to the arc spaces used to define the motivic zeta
 function.

We briefly sketch the ideas of the construction. Let $f$ be a
non-constant complex polynomial in $n$ variables, and let $x$ be a
closed point of the hypersurface $X_s$ defined by $f$. Let
$B=B(x,\varepsilon)$ be a small open disc around $x$ in $\C^n$,
and let $D=B(0,\eta)$ be a small open disc around the origin in
$\C$. We denote $D-0$ by $D^{\ast}$, and we denote by
$\widetilde{D^{\ast}}\rightarrow D^{\ast}$ the universal covering
space
$$\widetilde{D^{\ast}}=\{z\in\C\,|\,Im(z)>-\log (\eta)\}\rightarrow D^{\ast}:z\mapsto \exp(iz)$$
 For
$0<\eta\ll\varepsilon\ll 1$, the mapping
$$f:B\cap f^{-1}(D^\ast)\rightarrow D^{\ast}$$ is a locally trivial
fibration, called the Milnor fibration of $f$ at $x$. Its
canonical fiber $(B\cap
f^{-1}(D^\ast))\times_{D^{\ast}}\widetilde{D^{\ast}}$ is called
the canonical topological Milnor fiber $F_x$ of $f$ at $x$. The
group of deck transformations $\pi_1(D^{\ast})$ acts on $F_x$, and
the canonical generator $z\mapsto z+2\pi$ of $\pi_1(D^{\ast})$
induces a monodromy automorphism on the singular cohomology spaces
$H^i_{sing}(F_x,\C)$. The Milnor fibration is used as a tool to
gather information about the topology of $X_s$ near $x$ (see
\cite{Milnor}).

We can mimic this construction in the setting of formal geometry.
Let $k$ be an algebraically closed field of characteristic zero,
and let $f$ be a dominant morphism $f:X\rightarrow \mathrm{Spec}\,
k[t]$, with $X$ a smooth irreducible variety over $k$. Denote by
$X_s$ the hypersurface defined by $f$. We put $R=k[[t]]$, and
$K=k((t))$. Moreover, for each integer $d>0$, we denote by $K(d)$
the unique extension $k((\sqrt[d]{t}))$ of degree $d$ over $K$,
and we denote by $R(d)$ the normalization $k[[\sqrt[d]{t}]]$ of
$R$ in $K(d)$.

Taking the formal $t$-adic completion of $f$, we obtain a formal
$R$-scheme $\widehat{X}$, whose special fiber is isomorphic to
$X_s$. The ``complement'' of $X_s$ in $\widehat{X}$ is a smooth
rigid variety $X_\eta$ over $K$, and is endowed with a canonical
specialization map of ringed sites $sp:X_\eta\rightarrow
\widehat{X}$; see \cite[0.2]{bert}.

As we'll see in Section \ref{points}, for any integer $d>0$,
$K(d)$-valued points on $X_\eta$ correspond canonically to arcs
$\psi:\mathrm{Spec}\,R\rightarrow X$ on $X$, with $f(\psi)=t^d$.
 The specialization morphism $sp$ maps
$\psi$ to its origin $\psi(0)$. In Section \ref{milnor}, we define
the analytic Milnor fiber of $f$ at a closed point $x$ of $X_s$ as
$sp^{-1}(x)$, and we denote it by $\cF_x$. It is an open rigid
subspace of $X_\eta$. If $k=\C$, we prove in Section
\ref{cohomology} that its $\ell$-adic cohomology is isomorphic to
the singular cohomology of the topological Milnor fiber $F_x$, and
the absolute Galois action of $G(K^s/K)$ corresponds to the action
of the monodromy.

We define the local singular series $F(f;d)$ of $f$ by integrating
the Gelfand-Leray form associated to $f$ on $X_\eta\times_K K(d)$
(Definition \ref{localsing0} and \ref{localsing}; the definition
of the motivic integral of a volume form on a smooth rigid variety
is recalled in Subsection \ref{integral}). The result is a
function $F(f;.)$ from $\N^{*}$ to the localized Grothendieck ring
$\mathcal{M}_{X_s}$. Theorem \ref{weilzeta} states that the
motivic zeta function coincides (modulo a factor $\LL$) with the
Mellin transform $\sum_{d>0} F(f;d)T^d$.

Hence, we recover the motivic zeta function as a Weil-type zeta
function of the ``nearby fiber'' $X_\eta$. The corresponding
cohomological interpretation in Theorem \ref{dl} generalizes a
result by Denef and Loeser \cite[1.1]{DLLefschetz} which expresses
the Lefschetz numbers of the monodromy transformation, in terms of
the Euler characteristic of the jet spaces of the hypersurface
defined by $f$. We will obtain it as a special case of our Theorem
\ref{trace}, a Lefschetz trace formula for rigid spaces, which
relates the Euler characteristic of the motivic Serre invariant to
the Galois action on the \'etale cohomology.

Finally, we will prove in Section \ref{serreseries} that, for any
 generically smooth $stft$ formal $R$-scheme
$X_\infty$, and any gauge form $\omega$ on $X_\eta$, the volume
Poincar\'e series
$$S(X_\infty,\omega;T):=\left( \sum_{d>0}\int_{X_\eta\times_K K(d)}|\omega\otimes_K K(d)|\right) T^d\ \in\mathcal{M}_{X_s}[[T]]$$ is rational over $\mathcal{M}_{X_s}$,
 and we give an explicit expression.
This allows us to associate a motivic volume to
$X_\eta\widehat{\times}_K\widehat{K^s}$ in Section \ref{tame}.
When $X_\infty$ is the completion of a morphism $f:X\rightarrow
\A^1_k$ as above, this motivic volume coincides with Denef and
Loeser's motivic nearby cycles, as we will see in Section
\ref{motmilnor}.

To conclude this introduction, we give a survey of the structure
of this paper. In Section \ref{prel}, we recall some basic facts
on formal schemes and rigid varieties. Moreover, we prove
resolution of singularities for affine, flat, generically smooth
formal $k[[t]]$-schemes if $k$ has characteristic zero
(Proposition \ref{affineres}). In Section \ref{motrigsub}, we
recall the definition of the (relative) motivic Serre invariant,
as introduced in \cite{NiSe2}, and we briefly discuss motivic
Serre invariants with support. In Section \ref{weakner}, we
construct weak N\'eron models for the tame ramifications of a
regular flat $stft$ formal $R$-scheme with strict normal
crossings, where $R$ is a complete discrete valuation ring with
perfect residue field (Theorem \ref{neronmodel}). Our trace
formula for rigid spaces is proved in Section \ref{traceformula}
(Theorem \ref{trace}). In Section \ref{gauge}, we briefly recall
the definition of the motivic integral of a gauge form on a smooth
rigid variety. Next, we define the order of a gauge form $\omega$
on the generic fiber of a regular $stft$ formal $k[[t]]$-scheme
$\mX$ with strict normal crossings, assuming that $k$ has
characteristic zero (Definition \ref{order}), and we prove its
main properties. We compute the volume Poincar\'e series of the
pair $(\mX,\omega)$ in Section \ref{serreseries} (Theorem
\ref{expression}). The coefficients of this Poincar\'e series are
the motivic integrals of $\omega$ over the ramifications
$X_\eta\times_K K(d)$ (Definition \ref{volume}). The explicit
expression shows that the series is rational over the localized
Grothendieck ring $\mathcal{M}_{X_s}$, and allows us to define the
motivic volume of $X_\infty$ in Section \ref{tame} as a limit of
this Poincar\'e series (Definition \ref{motvolume}).

Finally, Section \ref{app} contains the applications to the theory
of motivic zeta functions. In Subsection \ref{milnor}, we define
the analytic Milnor fiber of a morphism of algebraic varieties
$f:X\rightarrow \A^1_k$, and we study its \'etale cohomology and
its points, in comparison to the the singular cohomology of the
topological Milnor fiber (if $k=\C$), and the points of the arc
space of $X$. In Subsection \ref{singseries}, we define the local
singular series of $f$, and we give an explicit expression in
terms of an embedded resolution of singularities for $f$. In
Subsection \ref{comparzeta}, we show that the Mellin transform of
the local singular series coincides with the motivic zeta function
of $f$ (Theorem \ref{weil}). To conclude, in Subsection
\ref{motmilnor}, we compare the motivic volume of the rigid nearby
fiber to Denef and Loeser's motivic nearby cycles.

The results in this paper were announced in \cite{NiSe-Milnor}. A
survey of the philosophy behind the project can be found in
\cite{NiSe3}.
\section{Preliminaries}\label{prel}
\subsection{Notation and conventions} Throughout this paper,
$R$ will be a complete discrete valuation ring, with quotient
field $K$, and perfect residue field $k$. We denote the
characteristic exponent of $k$ by $p$. We fix a uniformizing
element $t$ in $R$, i.e. a generator of the maximal ideal.

For any integer $d>0$, prime to $p$, we denote by $K(d)$ the
totally ramified extension $K[t_d]/((t_d)^d-t)$ of $K$, and by
$R(d)$ the normalization of $R$ in $K(d)$.
 For any $R$-scheme $X$, for any formal $R$-scheme $X_\infty$, and
for any integer $d>0$ prime to $p$, we write $X(d):=X\times_R
R(d)$ and $X(d)_\infty:=X_\infty\times_R R(d)$.

For any field $F$, we will denote by $F^s$ a separable closure. We
denote the normalization of $R$ in $K^s$ by $R^s$. The residue
field $k^s$ of $R^s$ is a separable closure of $k$. We'll write
$\widehat{K^s}$ for the completion of the valued field $K^s$, and
$\widehat{R^s}$ for the closure of $R^s$ in $\widehat{K^s}$.

Furthermore, we'll denote by $K^t$ the maximal tamely ramified
extension of $K$ in $K^s$, and by $R^t$ the normalization of $R$
in $K^t$. The residue field of $R^t$ is again equal to $k^s$.
We'll write $\widehat{K^t}$ for the completion of the valued field
$K^t$, and $\widehat{R^t}$ for the closure of $R^t$ in
$\widehat{K^t}$.

When $\cY$ is a rigid space over an algebraic, tamely ramified
extension $K'$ of $K$, we denote by $\overline{\cY}$ the base
change $\cY\widehat{\times}_{K'}\widehat{K^t}$.

A variety over a scheme $S$ is a reduced separated $S$-scheme of
finite type. We do not require it to be irreducible. For any
scheme $S$, we denote the underlying reduced scheme by $S_{red}$.
For any $k$-scheme $S$ of finite type, we denote by $Sm(S)$ its
$k$-smooth part.

For any $R$-scheme $X$, we denote its special fiber $X\times_R k$
by $X_s$. We call $X$ generically smooth, if its generic fiber
$X\times_R K$ is a smooth $K$-scheme.

If $T$ is any topological space, we write $H_{sing}(T,\C)$ for the
graded vector space $\oplus_i H^i_{sing}(T,\C)$, where
$H^i_{sing}(\,.\,)$ denotes the $i$-th singular cohomology space.

Recall the following convention: let $V=\oplus_{i\in \Z}V_i$ be a
graded vector space, with $V_i\neq 0$ for only a finite number of
$i$, and $V_i$ of finite dimension for all $i$. If $F=\oplus_i
F_i$ is a graded endomorphism of $V$, then the trace of $F$ on $V$
is defined as
$$Tr(F\,|\,V):=\sum_{i\in \Z}(-1)^i Tr(F_i\,|\,V_i)$$ Likewise,
the dimension of $V$ is the Euler characteristic
$$\mathrm{dim}\,V:=\sum_{i\in \Z}(-1)^i\mathrm{dim}\,V_i$$

We denote by $\N^{*}$ the set of integers $n>0$.
\subsection{Formal schemes and rigid varieties}
An $stft$ formal $R$-scheme $X_\infty$ is a separated formal
scheme, topologically of finite type over $R$. We denote its
special fiber by $X_s$, and its generic fiber (in the category of
separated quasi-compact rigid $K$-varieties) by $X_\eta$. We say
$X_\infty$ is generically smooth, if $X_\eta$ is a smooth rigid
$K$-variety. We denote by $Sm(X_\infty)$ the smooth part of
$X_\infty$ over $R$. If $R'$ is a finite extension of $R$, and
$\psi$ is an element of $X_\infty(R')$, then we denote by
$\psi(0)$ the image of the unique point of $\mathrm{Spf}\,R'$ in
$X_s$.

For any flat $R$-variety $X$, we denote its formal $t$-adic
completion by $\widehat{X}$. It is a $stft$ flat formal
$R$-scheme, and $\widehat{X}$ is generically smooth if $X$ is. If
$R=k[[t]]$, we will also write $\widehat{X}$ for the $t$-adic
completion of a morphism of $k$-varieties $X\rightarrow
\mathrm{Spec}\,k[t]$, flat over the origin.

As explained in \cite[0.2]{bert}, there exists a canonical
specialization morphism of ringed sites
$$sp:X_\eta\rightarrow X_\infty$$
For any locally closed subscheme $Z$ of $X_s$, we denote by
$\widehat{X_\infty/Z}$ the formal completion of $X_\infty$ along
$Z$. It is topologically of finite type over $R$ if $Z$ is open in
$X_s$, but not in general. Moreover, $sp^{-1}(Z)$ is a rigid
subspace of $X_\eta$, and this space is canonically isomorphic to
the generic fiber of $\widehat{X_\infty/Z}$. We denote this rigid
space by $]Z[$, and call it the tube of $Z$ in $X_\infty$. It is
quasi-compact if $Z$ is open in $X_s$, but not in general.

\subsection{Dilatations}
Let $X_\infty$ be a $stft$ flat formal $R$-scheme, and let $V$ be
a closed subscheme of $X_s$, defined by an ideal sheaf
$\mathcal{I}$ on $X_\infty$. We denote the formal blow-up of
$X_\infty$ at $V$ by $Bl_V:X'_\infty\rightarrow X_\infty$. The
dilatation $h:Y_\infty\rightarrow X_\infty$ of $X_\infty$ with
center $V$, is defined in \cite[\S 3]{formner} as the open formal
subscheme of $X'_\infty$ where $t$ generates
$\mathcal{I}\mathcal{O}_{X'_\infty}$. It has the following
universal property: $Y_\infty$ is $R$-flat, and
$h_s:Y_s\rightarrow X_s$ factors through $V\rightarrow X_s$.
Moreover, for any morphism of $stft$ formal $R$-schemes
$h':Z_\infty\rightarrow X_\infty$, such that $Z_\infty$ is
$R$-flat and $h'_s$ factors through $V$, there exists a unique
morphism $v:Z_\infty\rightarrow Y_\infty$ with $h'=h\circ v$. In
particular, for any finite unramified extension $R'$ of $R$, and
for any point $\psi$ in $X_{\infty}(R')$ with $\psi(0)\in V$, the
unique point of $X'_{\infty}(R')$ that maps to $\psi$ is contained
in $Y_{\infty}(R')$.

\subsection{Normal crossings and resolution of singularities}
\begin{definition} Let $X$ be a regular flat $R$-scheme (resp. a regular flat
$stft$ formal $R$-scheme). We say $X_s$ is a strict normal
crossing divisor, if we can find, for each closed point $x$ of
$X_s$, a regular system of local parameters $(x_0,\ldots,x_m)$ in
the local ring $\mathcal{O}_{X,x}$, such that
$t=u\prod_{i=0}^{m}x_i^{M_i}$, with $u$ a unit, and $M_i\in\N$.
\end{definition}

Let $X$ be a regular flat $R$-variety (resp. a regular flat $stft$
formal $R$-scheme) whose special fiber is a strict normal crossing
divisor. Let $E_i,\,i\in I$, be the irreducible components of
$(X_s)_{red}$. The variety $E_i$ is smooth over $k$, for each
$i\in I$. We'll denote by $N_i$ the length of the local ring of
$X_s$ at the generic point of $E_i$, and we call $N_i$ the
multiplicity of $E_i$ in $X_s$. We write $X_s=\sum_{i\in I}N_i
E_i$. We say that $X_s$ is a \textit{tame} strict normal crossing
divisor, if each $N_i$ is prime to the characteristic exponent $p$
of the residue field $k$.

For any non-empty subset $J$ of $I$, we define $E_J:=\cap_{i\in J}
E_i$ and $E_J^{o}:=E_J\setminus (\cup_{i\notin J}E_i)$. Moreover,
we put $m_J:=gcd\{N_i\,|\,i\in J\}$.

\begin{definition}
A resolution of singularities of a generically smooth flat
$R$-variety $X$ (resp. a generically smooth, flat $stft$ formal
$R$-scheme), is a proper birational morphism of flat $R$-varieties
(resp. a morphism of flat $stft$ formal $R$-schemes)
$h:X'\rightarrow X$, such that $h$ induces an isomorphism on the
generic fibers, and such that $X'$ is regular, with as special
fiber a strict normal crossing divisor $X'_s=\sum_{i\in I}N_i
E_i$. We say that the resolution $h$ is \textit{tame} if $X'_s$ is
a tame strict normal crossing divisor.
\end{definition}

\begin{lemma}\label{completion}
Let $A$ be an algebra topologically of finite type over $R$, and
let $X$ be a scheme of finite type over $R$. If $\mathfrak{M}$ is
a maximal ideal of $A$, defining a point $x$ of $\mathrm{Spf}\,A$,
then the local morphism $A_{\mathfrak{M}}\rightarrow
\mathcal{O}_{\mathrm{Spf}\,A,x}$ induces an isomorphism on the
completions (w.r.t. the respective maximal ideals)
$$\widehat{A}_{\mathfrak{M}}\cong\widehat{\mathcal{O}}_{\mathrm{Spf}\,A,x} $$
If $x$ is a point on the special fiber $X_s$ of $X$, then the
local morphism $\mathcal{O}_{X,x}\rightarrow
\mathcal{O}_{\widehat{X},x} $ induces a canonical isomorphism on
the completions
$$\widehat{\mathcal{O}}_{X,x}\cong\widehat{\mathcal{O}}_{\widehat{X},x} $$
\end{lemma}
\begin{proof}
The first point is shown in the proof of \cite[1.2.1]{conrad}. As
for the second, note that for any $n\geq 0$, and with
$X_n:=X\times_R (R/t^{n+1})$,
$$(\mathcal{O}_{X,x})/(t^{n+1})\cong \mathcal{O}_{X_n,x}\cong (\mathcal{O}_{\widehat{X},x})/(t^{n+1})$$
\end{proof}

\begin{lemma}\label{reg}
Let $A$ be an algebra topologically of finite type over $R$, and
let $X$ be a separated scheme of finite type over $R$.
\begin{enumerate}
\item $\mathrm{Spec}\,A$ is regular iff $\mathrm{Spf}\,A$ is
regular. Moreover, $X$ is regular iff $\widehat{X}$ is regular.
\item $(\mathrm{Spec}\,A)_s$ is a strict normal crossing divisor
iff $(\mathrm{Spf}\,A)_s$ is a strict normal crossing divisor.
Moreover, $X_s$ is a strict normal crossing divisor iff
$\widehat{X}_s$ is a strict normal crossing divisor. \item
$\mathrm{Spf}\,A$ is generically smooth, iff $A\otimes_R K$ is
geometrically regular over $K$. Moreover, if $X$ is generically
smooth, then $\widehat{X}$ is generically smooth. \item if $K$ is
perfect, any regular $stft$ formal $R$-scheme is generically
smooth.
\end{enumerate}
\end{lemma}
\begin{proof}
Regularity of a local Noetherian ring is equivalent to regularity
of its completion \cite[17.1.5]{ega4.1}, so (1) follows from Lemma
\ref{completion}. Point (2) follows from the fact that, for any
local Noetherian ring $S$, a tuple $(x_0,\ldots,x_m)$ in $S$ is a
regular system of local parameters for $S$ iff it is a regular
system of local parameters for $\widehat{S}$.

So let us prove (3).
Smoothness of the affinoid space
$\mathcal{Y}:=\mathrm{Sp}\,(A\otimes_R K)$ over $K$ is equivalent
to geometric regularity of $\mathcal{O}_{\mathcal{Y},x}$ over $K$
for each point $x$ of $\mathcal{Y}$, by \cite[2.8]{formrigIII}. If
we denote again by $x$ the corresponding closed point of
$Y:=\mathrm{Spec}\,(A\otimes_R K)$, then the natural local
morphism of local Noetherian rings $\mathcal{O}_{Y,x}\rightarrow
\mathcal{O}_{\mathcal{Y},x} $ induces an isomorphism on the
completions
$$\widehat{\mathcal{O}}_{\mathcal{Y},x}\cong
\widehat{\mathcal{O}}_{Y,x}$$ by
\cite[5.2.6/1,\,7.3.2/3,\,7.3.2/7]{BGR}.
 So the statement follows from the fact that
a local Noetherian ring is geometrically regular over $K$, iff its
completion is (by \cite[17.1.5]{ega4.1} and the fact that, for any
finite field extension $K'/K$ and any local Noetherian ring $C$
over $K$, $\widehat{C}\otimes_K K'\cong \widehat{C\otimes_K K'}$).

 If
$X$ is generically smooth, then so is $\widehat{X}$, since by
\cite[0.3.5]{bert}, the generic fiber $X_\eta$ of $\widehat{X}$ is
an open rigid subspace of the analytification of $X\times_R K$,
which is smooth by \cite[A.2.1]{conrad}.

Finally, to establish (4), it suffices to note that for any
affinoid algebra $B$, regularity of $B$ implies geometric
regularity of $B$ over $K$. This follows from
\cite[19.6.4-5]{ega4.1}.
\end{proof}

\begin{prop}\label{affineres}
If $k$ has characteristic zero, any affine generically smooth flat
$stft$ formal $R$-scheme $X_\infty=\mathrm{Spf}\,A$ admits a
resolution of singularities by means of formal admissible blow-ups
with smooth centers.
\end{prop}
\begin{proof}
First, we show that the scheme $X=\mathrm{Spec}\,A$ admits a
resolution of singularities $h:X'\rightarrow X$ by means of
blow-ups with smooth centers, concentrated in the special fiber
$X_s$. By \cite[8.2]{VMayor}, it suffices to prove that, for any
$n\geq 1$, the algebra $T_n:=R\{x_1,\ldots,x_n\}$ of convergent
power series over $R$ satisfies the following properties:
\begin{enumerate}
\item $Der_k(T_n)$ is a finite projective $(T_n)$-module, locally
of rank $(n+1)$, \item if $\mathfrak{M}$ is a maximal ideal of
$T_n$, then the dimension of the localization of $T_n$ at
$\mathfrak{M}$ equals $n+1$, and the residue field
$T_n/\mathfrak{M}$ is algebraic over $k$.
\end{enumerate}
Point (1) follows from the fact that $Der_k(T_n)\cong
(T_n)^{n+1}$, since a $k$-derivation on $T_n$ is determined by the
images of $t$ and the $x_i$ (such a derivation $D$ is
automatically continuous, since $D(t^i)$ is contained in the ideal
$(t^{i-1})$ by the Leibniz rule). Point (2) follows from the fact
that the completion of this localization is isomorphic to
$\widehat{\mathcal{O}}_{\A^n_R,x}$ for some closed point $x$ on
$\A^n_R$, by Lemma \ref{completion}.

By \cite[2.6.6]{bosch} and Lemma \ref{reg}, the $t$-adic
completion of $h$ is a resolution of singularities
$h:X'_\infty\rightarrow X_\infty$ by means of formal admissible
blow-ups with smooth centers.
\end{proof}
\subsection{\'Etale cohomology of rigid varieties}
Berkovich developed an \'etale cohomology theory for
non-archimedean analytic spaces (including the rigid spaces) in
\cite{Berk-etale}. Throughout this paper, all cohomology will be
Berkovich' \'etale cohomology, unless explicitly stated otherwise.
We fix a prime $\ell$, invertible in $k$. We define, for any
analytic space $\cY$ over the completion of some algebraic
extension of $K$,
\begin{eqnarray*}
H^{*}(\cY,\Z_{\ell})&=&\underleftarrow{\lim} H^{*}(\cY,\Z/\ell^n)
\\H^{*}(\cY,\Q_{\ell})&=&H^{*}(\cY,\Z_{\ell})\otimes\Q_{\ell}.
\end{eqnarray*}
We will simply write $H^{*}(\cY)$ for $H^{*}(\cY,\Q_{\ell})$, and
we'll denote the graded vector space $\oplus_i H^i(\cY)$ by
$H(\cY)$.

\section{Motivic Serre invariants with support}\label{motrigsub}
\subsection{The Grothendieck ring of varieties}
Let $Z$ be a variety over $k$. Consider the free Abelian group,
generated by the isomorphism classes $[X]$ of $Z$-varieties $X$.
We take the quotient of this group w.r.t. the following relations:
whenever $X$ is a $Z$-variety, and $Y$ is a closed subvariety of
$X$, we impose $[X]=[X\setminus Y]+[Y]$. This quotient is called
the Grothendieck group of varieties over $Z$, and is denoted by
$K_0(Var_Z)$. We denote the class $[\A^1_Z]$ of the affine line
over $Z$ by $\LL_Z$, or by $\LL$ if there is no risk of confusion.

A constructible subset $C$ of a $Z$-variety $X$ can be written as
a disjoint union of locally closed subsets, and defines
unambiguously an element $[C]$ of $K_0(Var_Z)$.
 When $Z$ is a separated scheme of finite
type over $k$, we will write $K_0(Var_Z)$ instead of
$K_0(Var_{Z_{\mathrm{red}}})$. For any separated scheme $X$ of
finite type over $Z$, we will write $[X]$ instead of
$[X_{\mathrm{red}}]$.

We can define a product on $K_0(Var_Z)$ as follows: for any pair
of $Z$-varieties $X$, $Y$, we put $[X].[Y]=[X\times_Z Y]$. This
definition extends bilinearly to a product on $K_0(Var_Z)$, and
makes it into a ring, the Grothendieck ring of varieties over $Z$.
The localized Grothendieck ring $\mathcal{M}_Z$ is obtained by
inverting $\LL_Z$ in $K_0(Var_Z)$.


A morphism of $k$-varieties $f:W\rightarrow Z$ induces
 base-change ring morphisms
$K_0(Var_Z)\rightarrow K_0(Var_W)$ and $\mathcal{M}_Z\rightarrow
\mathcal{M}_W$, as well as forgetful
 morphisms of Abelian groups $K_0(Var_W)\rightarrow K_0(Var_Z)$ and $\mathcal{M}_W\rightarrow
\mathcal{M}_Z$.

If $Z=\mathrm{Spec}\,k$, we write $K_0(Var_k)$, $\mathcal{M}_k$,
and $\LL$, rather than $K_0(Var_{\mathrm{Spec}\,k})$,
$\mathcal{M}_{\mathrm{Spec}\,k}$, and $\LL_{\mathrm{Spec}\,k}$.

The Grothendieck group $K_0(Var_Z)$ is a universal additive
invariant for $Z$-varieties: if $A$ is an Abelian group, and
$\chi$ is an invariant of $Z$-varieties taking values in $A$, such
that, for any $Z$-variety $X$ and any closed subvariety $Y\subset
X$, $\chi(X)=\chi(X\setminus Y)+\chi(Y)$, then $\chi$ factors
uniquely through a group morphism $\chi:K_0(Var_Z)\rightarrow A$,
defined by $\chi([X])=\chi(X)$. If $A$ is a ring, and $\chi$ is
multiplicative, i.e. $\chi((X\times_Z
Y)_{\mathrm{red}})=\chi(X).\chi(Y)$ for any pair of $Z$-varieties
$X,Y$, then $\chi:K_0(Var_Z)\rightarrow A$ is a morphism of rings.

For instance, for any $k$-variety $X$, we can consider its
topological Euler characteristic $\chi_{top}(X)$. Fix a prime
$\ell$, invertible in $k$. Then $\chi_{top}(X)$ is defined as
$$\chi_{top}(X):=\sum_{i\geq 0}(-1)^i\,\mathrm{dim}\,H_c^i(X\times_k k^s,\Q_\ell) $$
where $H^i_c(\, .\, ,\Q_\ell)$ is $\ell$-adic \'etale cohomology
with proper support, and $k^s$ is a separable closure of $k$. This
is an additive invariant, hence defines a morphism of groups
$$\chi_{top}:K_0(Var_Z)\rightarrow \Z $$ for any $k$-variety $Z$.
It is multiplicative for $Z=\mathrm{Spec}\,k$, so we get a
morphism of rings
$$\chi_{top}:K_0(Var_k)\rightarrow \Z .$$

\subsection{Motivic Serre invariants with support}
 Let $X_\eta$
be a separated, quasi-compact smooth rigid $K$-variety. A
\textit{weak N\'eron $R$-model for $X_\eta$} is a smooth $stft$
formal scheme $U_\infty$ over $Spf\,R$, whose generic fiber is an
open rigid
 subspace of $X_\eta$, and which has the property that the natural
map $U_\infty(R^{sh})\rightarrow X_\eta(K^{sh})$ is bijective
\cite[Definition 1.3]{formner}. By this latter property, we mean
that $U_\infty(R')\rightarrow X_\eta(K')$ is bijective, for any
finite unramified extension $K'$ of $K$, where $R'$ is the
normalization of $R$ in $K'$. Observe that this map is always
injective, since $U_\infty$ is separated.

\begin{definition}
Let $h:U_\infty\rightarrow X_\infty$ be a morphism of $stft$ flat
formal $R$-schemes, with $X_\eta$ smooth over $K$. We say that $h$
is a N\'eron $R$-smoothening for $X_\infty$, if it has the
following properties:
\begin{enumerate} \item $U_\infty$ is a weak N\'eron $R$-model for
$X_\eta$, \item there exists a morphism
 $X'_\infty\rightarrow X_\infty$  of $stft$ flat formal $R$-schemes,
 inducing an isomorphism $X'_\eta\rightarrow X_\eta$ on the generic
 fibers, such that
 $h$ factors through
an open immersion $U_\infty\rightarrow X'_\infty$.
\end{enumerate} 
\end{definition}

The result in \cite[\S3, Theorem 3.1]{formner}
can be interpreted in our context as follows~:

\begin{theorem}\label{bs}
If $X_\infty$ is a generically smooth, $stft$ flat formal
$R$-scheme, then $X_\infty$ admits a N\'eron $R$-smoothening.
Moreover, we can always find an admissible blow-up
$X'_\infty\rightarrow X_\infty$, such that
$Sm(X'_\infty)\rightarrow X_\infty$ is a N\'eron $R$-smoothening.
\end{theorem}

In \cite{NiSe2}, we refined the notion of motivic Serre invariant,
first introduced in \cite{motrigid}, as follows:

\begin{definition}[Motivic Serre invariant]\label{motserre}
If $X_\infty$ is a generically smooth, $stft$ flat formal
$R$-scheme, the \textit{(relative) motivic Serre invariant}
$S(X_\infty)$ of $X_\infty$, is the class $[U_s]$ in
$K_0(Var_{X_s})/(\LL_{X_s}-[X_s])$, where $U_\infty\rightarrow
X_\infty$ is any N\'eron $R$-smoothening of $X_\infty$.
\end{definition}

We proved in \cite{NiSe2}, Theorem 6.1, that this definition does
not depend on the choice of the N\'eron $R$-smoothening, in the
case where $X_\eta$ has pure dimension. We proved the general case
in \cite{NiSe-weilres}, Theorem 5.9.

If $Y_\infty \rightarrow X_\infty$ is a morphism of generically
smooth, $stft$ flat formal $R$-schemes, inducing an isomorphism on
the generic fibers,  the forgetful morphism
$$K_0(Var_{Y_s})/(\LL_{Y_s}-[Y_s])\rightarrow K_0(Var_{X_s})/(\LL_{X_s}-[X_s])$$ maps
$S(Y_\infty)$ to $S(X_\infty)$.  Hence, the motivic Serre
invariants computed over all $stft$ flat formal $R$-models of a
separated, quasi-compact, smooth rigid variety $X_\eta$ over $K$,
form a projective system.

For any $stft$ flat formal $R$-model $X_\infty$ of $X_\eta$, the
Serre invariant $S(X_\eta)$ defined in \cite{motrigid}, is the
image of $S(X_\infty)$ under the forgetful morphism
$$K_0(Var_{X_s})/(\LL_{X_s}-[X_s])\rightarrow
K_0(Var_{k})/(\LL-1)$$ It only depends on the rigid space
$X_\eta$, and not on the choice of the formal model $X_\infty$.

\begin{definition}
For any locally closed subscheme $V$ of $X_s$, we define the
motivic Serre invariant $S_V(X_\infty)$ of $X_\infty$ with support
in $V$ as the image of $S(X_\infty)$ under the base change
morphism $$K_0(Var_{X_s})/(\LL_{X_s}-[X_s]) \rightarrow
K_0(Var_{V})/(\LL_{V}-[V])$$
\end{definition}

We will think of $S_{V}(\mX)$ as a measure for the number of
unramified points on the tube $]V[$.

In general, this tube is not quasi-compact, so we cannot take its
motivic Serre invariant in a direct way.

\begin{prop}

\begin{enumerate}
\item If $h:U_\infty\rightarrow X_\infty$ is a N\'eron
$R$-smoothening, then $$S_V(X_\infty)=[U_s\times_{X_s}V]\in
K_0(Var_{V})/(\LL_{V}-[V])$$
 \item If $Y_\infty$ is any open formal subscheme of $X_\infty$ containing $V$, then $S_V(\mX)=S_V(Y_\infty)$.
In particular, if $V$ is open in $X_s$, then
$S_V(\mX)=S(\widehat{\mX/V})$. \item Suppose that $V$ is closed in
$X_s$, and denote by $\pi:Y_\infty\rightarrow X_\infty$ the
dilatation with center $V$. Then $S_V(\mX)$ is the image of
$S(Y_\infty)$ under the forgetful morphism
$$K_0(Var_{Y_s})/(\LL_{Y_s}-[Y_s])\rightarrow  K_0(Var_{V})/(\LL_{V}-[V])$$
\end{enumerate}
\end{prop}
\begin{proof}
Point $(1)$ is clear. Point $(2)$ follows from the fact that, if
$h:U_\infty\rightarrow X_\infty$ is a N\'eron $R$-smoothening,
then $h:h^{-1}(Y_\infty)\rightarrow Y_\infty$ is one, as well.

Finally, as for $(3)$, let $Bl_V: X'_\infty\rightarrow X_\infty$
be the formal blow-up of of $\mX$ at $V$. The dilatation
$Y_\infty$ is an open formal subscheme of $X'_\infty$. Take a
N\'eron $R$-smoothening $h:U_\infty\rightarrow X'_\infty$. If we
denote the complement of $Y_s$ in $Bl_V^{-1}(V)$ by $E$, it
suffices to prove that $(U_\infty\setminus h^{-1}(E))\rightarrow
X'_\infty$ is still a N\'eron $R$-smoothening. This follows from
the fact that the tube $]E[$ in $X'_\infty$ does not contain any
$K'$-points, with $K'$ finite and unramified over $K$, by the
universal property of the dilatation. Computing $S_V(\mX)$ on this
N\'eron smoothening, we get
$$S_V(\mX)=[h^{-1}(Y_s)] =S(\mY)$$
in $K_0(Var_{V})/(\LL_{V}-[V])$.
\end{proof}

\section{Weak N\'eron models for ramifications}\label{weakner}
Let $\mX$ be a regular flat $stft$ formal $R$-scheme, such that
its special fiber $X_s$ is a strict normal crossing divisor
$\sum_{i\in I}N_i E_i$. We fix an integer $d>0$, prime to $p$. The
aim of this section, is to construct a N\'eron $R$-smoothening for
the
 ramification $X(d)_\infty:= \mX\times_R R(d)$.

Choose any subset $J$ of $I$, such that $m_J$ is prime to $p$. We
can cover $E_J^o$ by affine open formal subschemes
$U_\infty=\mathrm{Spf}\,V$, such that on $U_\infty$, we can write
$t=u\prod_{i\in J}x_i^{N_i}$, with $u$ a unit. We define an
\'etale cover $U_\infty'\rightarrow U_\infty$ by
$U_\infty'=\mathrm{Spf}\,V\{T\}/(uT^{m_J}-1)$. These covers glue
together, and we obtain an \'etale cover
$\widetilde{E}_J^o\rightarrow E_J^o$.

\begin{remark}\label{compare}
If $k$ has characteristic zero, and $\mX$ is isomorphic to the
formal completion of $Z\times_{k[t]} k[[t]]$, for some smooth
irreducible $k$-variety $Z$, and some dominant morphism
$Z\rightarrow \mathrm{Spec}\,k[t]$, then our cover
$\widetilde{E}_J^o/E_J^o$ coincides with the one defined in
\cite[2.3]{DLLefschetz}.
\end{remark}

We denote by $\widetilde{X(d)_\infty}\rightarrow X(d)_\infty$ the
normalization of $X(d)_\infty$ (see \cite{conrad} on normalization
of formal schemes), and we denote by $\widetilde{E(d)}_i^o$ the
inverse image $\widetilde{X(d)}_s\times_{X_s}E_i^o $ of $E_i^o$ in
$\widetilde{X(d)}_s$, for each $i\in I$.

\begin{definition}\label{elinear}
Let $J$ be a non-empty subset of $I$. We say that an integer
$d\geq 1$ is $J$-\textit{linear} if there exists, for each $j\in
J$, an integer $\alpha_j\in\N^{\ast}$, such that
$$
d=\sum_{j\in J}\alpha_j N_{j}
$$
 We say that an integer $d\geq 1$ is $X_s$-\textit{linear}
if there exists a subset $J\subset I$, with $\vert J\vert>1$ and
$E_J^o\neq \emptyset$, such that $d$ is $J$-linear.
\end{definition}

We recall two Lemmas from \cite{NiSe2}.
\begin{lemma}[\cite{NiSe2}, Lemma 5.15]\label{valuation}
Let $R'$ be a finite extension of $R$ of ramification degree $d$,
and consider an element $\psi$ of $X_\infty(R')$. If $J$ is the
unique subset of $I$ with $\psi(0)\in E_J^{o}$, then $d$ is
$J$-linear.
\end{lemma}
\begin{lemma}[\cite{NiSe2}, Lemma 5.17]\label{linear}
There exists a sequence of admissible blow-ups
$\pi^{(j)}:\mX^{(j+1)}\rightarrow \mX^{(j)}$, $j=0,\ldots,r-1$,
such that
\begin{itemize}
\item $\mX^{(0)}=\mX$, \item the special fiber of $\mX^{(j)}$ is a
strict normal crossing divisor
$$X_s^{(j)}=\sum_{i\in I^{(j)}} N_i^{(j)}E_i^{(j)},$$ \item $\pi^{(j)}$
is the formal blow-up with center $E^{(j)}_{J^{(j)}}$, for some
 subset $J^{(j)}$ of $I^{(j)}$, with $|J^{(j)}|>1$, \item $d$ is not
$X_s^{(r)}$-linear.
\end{itemize}
\end{lemma}

\begin{lemma}\label{normal}
 If
$d$ is not $X_s$-linear, then
$$Sm(\widetilde{X(d)}_s)=\bigsqcup_{N_i|d}\widetilde{E(d)}_i^o$$
Moreover, for any $i\in I$ with $N_i | d$, the morphism
$\widetilde{X(d)_\infty}\rightarrow \widetilde{X(N_i)_\infty}$
induces an isomorphism $$\widetilde{E(d)}_i^o\cong
\widetilde{E(N_i)}_i^o$$ and the $E_i^o$-variety
$\widetilde{E(N_i)}_i^o$ is canonically isomorphic to
$\widetilde{E}_i^o$.
\end{lemma}
\begin{proof}
Let $x$ be a smooth point of $\widetilde{X(d)}_s$. Since
$\widetilde{X(d)_\infty}$ is flat, this implies that
$\widetilde{X(d)_\infty}\rightarrow \mathrm{Spf}\,R(d)$ is smooth
at $x$, and hence, there exists a section $\psi:\mathrm{Spf}\,
R'\rightarrow \widetilde{X(d)_\infty}$, with $R'$ an unramified
extension of $R(d)$, and $\psi(0)=x$. If we denote by $J$ the
unique subset of $I$ such that the image of $x$ in $X_s$ belongs
to $E_J^o$, then Lemma \ref{valuation}, combined with the
hypothesis that $d$ is not $X_s$-linear, implies that $J$ is a
singleton $\{i\}$, and $N_i|d$. Hence,
$$Sm(\widetilde{X(d)}_s)\subset \bigsqcup_{N_i|d}\widetilde{E(d)}_i^o$$

Conversely, it follows from \cite{Kempf} (proof of the semi-stable
reduction theorem II, pages 198--202) that
$$\bigsqcup_{N_i|d}\widetilde{E(d)}_i^o \subset Sm(\widetilde{X(d)}_s)$$
(their arguments in the algebraic setting carry over to formal
schemes).

In fact, we can give an explicit description of
$Sm(\widetilde{X(d)_\infty})$.
 Choose $i\in I$
such that $N_i|d$, and let $x$ be any closed point of $E_i^o$.
Choose an affine open formal neighborhood
$U_\infty=\mathrm{Spf}\,V$ of $x$ in $X_\infty$, such that on
$U_\infty$, we can write $t=ux_i^{N_i}$, with $u$ a unit. For any
minimal prime ideal $\mathfrak{P}$ of $V$, we can write
$$
\left(\frac{x_i}{(t_d)^{d/N_i}}\right)^{N_i} -u^{-1}=0
$$ in the quotient field of $V/\mathfrak{P}$,
 so the normalization map $\widetilde{U(d)_\infty}\rightarrow
U(d)_\infty$ factors through $$\widetilde{U(d)_\infty}\rightarrow
Y(d):=\mathrm{Spf}\,V(d)\{T\}/((t_d)^{d/N_i}T-x_i,
uT^{N_i}-1)\rightarrow U(d)_\infty
$$ where $V(d):=V\otimes_R R(d)$.
We'll show that $Y(d)$ is smooth over $R(d)$. In particular,
$Y(d)$ is normal, so $\widetilde{U(d)_\infty}\rightarrow Y(d)$ is
an isomorphism.

Since $Y(d)$ is flat over $R(d)$, it suffices to show that the
special fiber $Y(d)_s$ is smooth over $k$. Reduction modulo $t_d$
yields $$Y(d)_s=\mathrm{Spec}\,V[T]/(x_i,uT^{N_i}-1)$$ The section
$x_i$ is part of a regular system of local parameters, so the
scheme $\mathrm{Spec}\,V/(x_i)$ is regular. Since $k$ is perfect,
$\mathrm{Spec}\,V/(x_i)$ is smooth over $k$. But $Y(d)_s$ is
\'etale over $\mathrm{Spec}\,V/(x_i)$, and hence smooth over $k$.

Finally, the explicit description $\widetilde{U(d)_\infty}\cong
\mathrm{Spf}\,V(d)\{T\}/((t_d)^{d/N_i}T-x_i,uT^{N_i}-1)$ shows
that $\widetilde{E(d)}^o_i$ does not depend on $d$, as long as
$N_i|d$. In fact, the restriction of $\widetilde{E(d)}^o_i$ over
$U_s$ is given explicitly by
$\mathrm{Spec}\,V[T]/(x_i,uT^{N_i}-1)$, which is canonically
isomorphic to the restriction of $\widetilde{E}_i^o$ over $U_s$.
\end{proof}

\begin{theorem}\label{neronmodel}
Let $\mX$ be a regular, generically smooth, flat $stft$ formal
$R$-scheme, such that its special fiber $X_s$ is a strict normal
crossing divisor $\sum_{i\in I}N_i E_i$. Let $d>0$ be an integer,
prime to $p$, such that $d$ is not $X_s$-linear. Then
$$Sm(\widetilde{X(d)_\infty})\rightarrow X(d)_\infty$$ is an
N\'eron $R$-smoothening, and
$$S(X(d)_\infty)=\sum_{N_i|d}[\widetilde{E}_i^o]$$
in $K_0(Var_{X_s})/(\LL_{X_s}-[X_s])$.
\end{theorem}
\begin{proof}
We only have to prove that $Sm(\widetilde{X(d)_\infty})$ is a weak
N\'eron model for $X(d)_\eta$. This follows from Lemma
\ref{valuation}, and Lemma \ref{normal}.
\end{proof}
\section{A trace formula for non-archimedean analytic
spaces}\label{traceformula} The purpose of this section is to
prove a Grothendieck trace formula for non-archimedean analytic
spaces.

We suppose that the residue field $k$ is algebraically closed. Let
$\varphi$ be a topological generator of the tame Galois group
$G(K^t/K)$.

When $\cY$ is a rigid space over an algebraic, tamely ramified
extension $K'$ of $K$, we denote by $\overline{\cY}$ the base
change $\cY\widehat{\times}_{K'}\widehat{K^t}$. Fix a prime
$\ell$, different from $p$.

%
%
%
\begin{lemma}\label{illusie}
 Let $F$ be a separably
closed field, and let $G=<g>$ be a finite cyclic group.
Let $Y$ be a normal proper $F$-variety, and let $U$ be an open
subscheme of $Y$. Let $\mathcal{L}$ be a lisse constructible
$\overline{\mathbb{Q}}_{\ell}$-sheaf on $U$, tamely ramified on
$Y-U$. Suppose that $G$ acts on $\mathcal{L}$.
 For
any closed point $x$ on $U$,
$$Tr(\,g\,|\,H_c(U,\mathcal{L}))=\chi_{top}(U)Tr(\,g\,|\,\mathcal{L}_x).$$
\end{lemma}
\begin{proof}
By definition, $\mathcal{L}$ is a lisse constructible sheaf over
some finite extension $Q$ of $\Q_{\ell}$.
 The tamely ramified lisse sheaf
$\mathcal{L}$ on $U$ with stalk $\mathcal{L}_x$ at $x$ is
determined by a continuous morphism
$$\psi:\pi_1^t(U,x)\rightarrow Aut\,_{Q}(\mathcal{L}_x).$$
Giving an endomorphism $f$ of $\mathcal{L}$, amounts to giving an
element $\rho(f)$ of $End\,_{Q}(\mathcal{L}_x)$, such that
$\rho(f)\psi(\gamma)=\psi(\gamma)\rho(f)$, for each $\gamma\in
\pi_1^t(U,x)$.

Passing to a finite extension of $Q$, we can decompose
$\mathcal{L}$ as $\oplus_i\mathcal{L}_i$, according to the
decomposition of $\mathcal{L}_x\otimes \overline{\Q}_{\ell}$ into
generalized eigenspaces w.r.t. $\rho(g)$. Hence, we may as well
assume that $\rho(g)$ has only one eigenvalue $e$ on
$\mathcal{L}_x\otimes \overline{\Q}_{\ell}$.

Write the automorphism induced by $\rho(g)$ on $\mathcal{L}_x$ as
$e.Id+N$, with $N$ nilpotent. Since
$Tr(\,N\,|\,H_c(U,\mathcal{L}))=0$, it suffices to prove that
$$e.\chi(U,\mathcal{L})=Tr(e.Id\,|\,H_c(U,\mathcal{L}))=\chi(U,Q)Tr(e.Id\,|\,\mathcal{L}_x)
=\chi_{top}(U).e.rank(\mathcal{L}_x).$$ This follows from
\cite{Illusie}, Cor 2.7.
\end{proof}

\begin{lemma}\label{zero}
Let $Y$ be a regular flat variety over $R$, such that $Y_s$ is a
strict normal crossing divisor $\sum_{i\in I} N_i E_i$. If $J$ is
a subset of $I$, with $|J|>1$, and $Bl_{E_J}$ is the blow-up of
$Y$ at $E_J$, with exceptional component $E'_0$, then $(E'_0)^o$
is a Zariski-locally trivial fibration over $E^o_J$, and its fiber
is a torus $\G_{m,k}^{|J|-1}$. In particular,
$[(E'_0)^o]=[E^o_J](\LL_{X_s}-[X_s])^{|J|-1}$ in $K_0(Var_{X_s})$.
\end{lemma}
\begin{proof}
 Let $x$ be any closed point
on $E_J$, and take a regular system of local parameters
$(x_0,\ldots,x_m)$ on $Y$ at $x$, such that $t=u\prod_{i\in
J}x_i^{N_i}$, with $u$ a unit, where we identified $J$ with a
subset of $\{0,\ldots,m\}$. The sequence $(x_i)_{i\in J}$ is
regular in a neighborhood $U$ of $x$. Now apply
\cite[IV-26]{Eisenbud},
 to the blow-up of $U$ with center $(x_i)_{i\in J}$.
\end{proof}
\begin{lemma}\label{nearby}
Let $Y$ be a regular flat variety over $R$, such that $Y_s$ is a
tame strict normal crossing divisor $\sum_{i\in I} N_i E_i$. The
complex of $\ell$-adic tame nearby cycles $R\psi^t_\eta(\Q_\ell)$
is constructible and tamely ramified on $Y_s$, and lisse on the
strata $E_J^o$. For any closed point $y$ of $Y_s$, and for any
integer $d>0$ prime to $p$,
\begin{eqnarray*}
Tr(\varphi^d\,|\,R\psi^t_\eta(\Q_\ell)_y)&=& N_i\quad \mbox{if}\
y\in E_i^o\ \mbox{with}\ N_i|d,
\\ Tr(\varphi^d\,|\,R\psi^t_\eta(\Q_\ell)_y)&=& 0\quad \mbox{else.}
\end{eqnarray*}
\end{lemma}
\begin{proof}
This follows from the explicit computation of the nearby cycles in
\cite[I,3.3]{sga7a}.
\end{proof}

\begin{theorem}[Trace formula]\label{trace}
Let $X$ be a generically smooth flat variety over $R$, and suppose
that $X$ admits a tame resolution $h:Y\rightarrow X$, with
$Y_s=\sum_{i\in I} N_i E_i$. Let $Z$ be a proper subvariety of
$X_s$. For any integer $d>0$, prime to $p$, we have
$$\chi_{top}\left(S_{Z}(\widehat{X(d)})\right)
=Tr(\varphi^d\,|\,H(\,\overline{]Z[}\,))=\sum_{N_i |d} N_i
\chi_{top}(h^{-1}(Z)\cap E_i^o).$$
\end{theorem}
\begin{proof}
We may suppose that $X=Y$.

By Berkovich' quasi-isomorphism \cite[3.5]{berk-vanish2},
$$Tr(\varphi^d\,|\,H(\,\overline{]Z[}\,))=Tr(\varphi^d\,|\,H(Z,R\psi^t_{\eta}(\Q_{\ell}))).$$
By Lemma \ref{illusie}, and Lemma \ref{nearby},
$$Tr(\varphi^d\,|\,H(\,\overline{]Z[}\,))=\sum_{N_i |d} N_i
\chi_{top}(Z\cap E_i^o),$$.

We use Lemma \ref{linear} to construct a resolution $Y'\rightarrow
X$, with $Y'_s=\sum_{i\in I'} N'_i E'_i$, such that $d$ is not
$Y'_s$-linear. Denote by $Z'$ the inverse image of $Z$ in $Y'_s$.

 By Theorem \ref{neronmodel},
$$\chi_{top}(S_{Z}(\widehat{X(d)}))
=\sum_{N'_i |d} \chi_{top}(\widetilde{Z}_i),$$ where
$\widetilde{Z}_i$ denotes the inverse image of $Z'\cap (E'_i)^o$
in $\widetilde{E'_i}^o$. Since $\widetilde{E'_i}^o$ is \'etale
over $(E'_i)^o$, of degree $N'_i$, and tamely ramified,
$\chi_{top}(\widetilde{Z}_i)=N'_i\chi_{top}(Z'\cap (E'_i)^o)$, for
any $i$. When $E'_i$ is an exceptional component of $Y'\rightarrow
Y$, then $\chi_{top}(Z'\cap (E'_i)^o)=0$, by Lemma \ref{zero}.
Hence,
$$\chi_{top}(S_{Z}(\widehat{X(d)}))
=\sum_{N_i |d} N_i \chi_{top}(Z\cap E_i^o).$$
\end{proof}

The conditions of Theorem \ref{trace} are satisfied, in
particular, when $k$ has characteristic zero.

\begin{cor}
Suppose that $k$ is an algebraically closed field of
characteristic zero. If $X^{an}_K$ is the analytification of a
smooth, proper $K$-variety $X_K$, then
$$\chi_{top}\left(S(X_K^{an}\times_K K(d))\right)=Tr(\,\varphi^d\,|\,H(X_K\times_K K^s,\Q_\ell))$$
for any integer $d>0$.
\end{cor}
\begin{proof}
If $X$ is any flat $R$-model for $X_K$, then $X_K^{an}$ is
canonically isomorphic to the generic fiber $X_\eta$ of
$\widehat{X}$, by \cite[0.3.5]{bert}. Moreover, by
\cite[7.5.4]{Berk-etale}, there is a canonical isomorphism
$$H^i(X_K\times_K K^s,\Q_\ell)\cong H^i(\overline{X_\eta})$$ for each $i\geq 0$.
\end{proof}
\begin{remark}
That some tameness condition is needed in the statement of the
Trace Formula, is already clear from the following example: let
$R$ be the ring $W(\mathbb{F}_p^s)$ of Witt vectors over
$\mathbb{F}_p^s$, where $p>0$ is a prime. Let $X_\infty$ be the
formal $R$-scheme $Spf\,R\{x\}/(x^p-p)$. Obviously, $S(X_\eta)=0$,
while
$Tr(\varphi\,|\,H(X_\eta\widehat{\times}_K\widehat{K^{t}}))=1$,
when $\varphi$ is a topological generator of the tame Galois group
$G(K^t/K)$.

It should be possible to replace $X$ in the statement of the trace
formula by an arbitrary generically smooth $stft$ formal
$R$-scheme $X_\infty$ which admits a tame resolution of
singularities. Also the condition that $Z$ is proper does not seem
essential. In particular, we expect that
$$\chi_{top}\left(S(X_\eta\times_K K(d))\right)=Tr(\,\varphi^d\,|\,H(\overline{X_\eta}))$$
holds for any separated smooth quasi-compact rigid variety
$X_\eta$ over $K$, if $k$ has characteristic zero (and also for
any separated smooth rigid variety $X_\eta$ over $K$ which can be
realized as a tube in a $stft$ formal $R$-scheme, using the
motivic Serre invariant with support in the left hand side).
\end{remark}

\section{Order of a gauge form}\label{gauge}
We assume $R=k[[t]]$, with $k$ a field of characteristic zero. Let
$X_\infty$ be a regular, flat $stft$ formal $R$-scheme, of pure
relative dimension $m$, such that the special fiber $X_s$ is a
strict normal crossing divisor $\sum_{i\in I}N_i E_i$. Let
$\omega$ be a gauge form on the generic fiber $X_\eta$. The
purpose of this section, is to define the order $ord_{E_i}\omega$
of $\omega$ at the generic point of $E_i$, for any $i\in I$.

\subsection{Motivic integral of a gauge form on a smooth rigid
variety}\label{integral} If $\mX$ is smooth, this order was
defined already in \cite[4.3]{formner}, for general $R$. This
definition was used in \cite[4.3.1]{motrigid} to give an
expression for the integral of a gauge form on a separated,
smooth, quasi-compact rigid variety over $K$. In \cite{NiSe2}, we
refined this notion as follows. Let $\mY$ be a generically smooth,
flat $stft$ formal $R$-scheme, of pure relative dimension $m$, and
let $\omega$ be a gauge form on $Y_\eta$. We take a N\'eron
$R$-smoothening $h:\mZ\rightarrow \mY$, and we denote by
$\mathcal{C}$ the set of connected components of $Z_s$. Then
$$\int_{Y_\infty}|\omega|:=\LL^{-m}\sum_{U\in
\mathcal{C}}[U]\LL^{-ord_U(h^*\omega)}\in \mathcal{M}_{Y_s}$$
depends only on the pair $(\mY,\omega)$, and not on $\mZ$ (see
\cite{NiSe2}, Lemma 6.4). The image of $\int_{Y_\infty}|\omega|$
under the forgetful morphism $\mathcal{M}_{Y_s}\rightarrow
\mathcal{M}_k$ is the motivic integral $\int_{Y_\eta}\omega$ from
\cite[4.1.2]{motrigid}. It depends only on the pair
$(Y_\eta,\omega)$, and not on the model $\mY$. Finally, the image
of $\int_{Y_\infty}|\omega|$ under the projection morphism
$$\mathcal{M}_{Y_s}\rightarrow
\mathcal{M}_{Y_s}/(\LL_{Y_s}-[Y_s])\cong
K_0(Var_{Y_s})/(\LL_{Y_s}-[Y_s])$$ is exactly the motivic Serre
invariant $S(\mY)$ from Definition \ref{motserre}. In particular,
it depends only on $\mY$, and not on $\omega$.
\subsection{The order of a top form at a section}
First, we generalize a definition from \cite[4.1]{motrigid}. Let
$Y_\infty$ be any flat $stft$ formal $R$-scheme, equidimensional
of relative dimension $m$. Let $R'$ be a finite extension of $R$,
of ramification index $e$.

\begin{definition}\label{orderideal}
For any element $\psi$ of $Y_\infty(R')$, and any ideal sheaf
$\mathcal{I}$ on $Y_\infty$, we define $ord(\mathcal{I})(\psi)$ as
the length of the $R'$-module $R'/\psi^*\mathcal{I}$.
\end{definition}

We recall that the length of the zero module is $0$, and the
length of $R'$ is $\infty$.

For any element $\psi$ of $Y_\infty(R')$, the $R'$-module
$M:=(\psi^{*}\Omega^m_{Y_\infty/R})/(\mathrm{torsion})$ is free of
rank $1$.

\begin{definition}
For any global section $\omega$ of $\Omega^m_{Y_\infty/R}$, we
define the order of $\omega$ at $\psi$ as the length of the
$R'$-module $M/R'(\psi^*\omega)$. We denote this value by
$ord(\omega)(\psi)$.
\end{definition}

If $e=1$, this definition coincides with the one given in
\cite[4.1]{motrigid}. It only depends on an open formal
neighbourhood of $\psi(0)$ in $Y_\infty$.

If $\omega\in \Omega^m_{Y_\eta/K}(Y_\eta)$, there exists an
integer $a\geq 0$ such that $t^a\omega\in
\Omega^m_{Y_\infty/R}(Y_\infty)$, by the isomorphism of sheaves
\cite[1.5]{formrigIII}
$$\Omega^m_{Y_\eta/K}\cong \Omega^m_{Y_\infty/R}\otimes_R K$$ and the fact
that $Y_\infty$ is quasi-compact.

\begin{definition}
 If $\omega$ is a global section of
$\Omega^m_{Y_\eta/K}$, we take an integer $a\geq 0$ such that
$t^a\omega$ is defined on $Y_\infty$, and we define
$ord(\omega)(\psi)$ as $ord(t^a\omega)(\psi)-ea$.
\end{definition}
This definition does not depend on the choice of $a$. If $Y_\eta$
is smooth, and $\omega$ is a gauge form on $Y_\eta$,
$ord(\omega)(\psi)$ is finite.


Now let $h:Z_\infty\rightarrow Y_\infty$ be a morphism of flat,
generically smooth $stft$ formal $R$-schemes, both equidimensional
of relative dimension $m$. Let $R'$ be a finite extension of $R$,
and fix a section $\psi$ in $Z_\infty(R')$. The canonical morphism
$$h^*\Omega^m_{Y_\infty/R}\rightarrow \Omega^m_{Z_\infty/R}$$
induces a morphism of free rank $1$ $R'$-modules
$$(\psi^*h^*\Omega^m_{Y_\infty/R})/(\mathrm{torsion})\rightarrow (\psi^*\Omega^m_{Z_\infty/R})/(\mathrm{torsion})$$
We define $ord(Jac_h)(\psi)$ as the length of its cokernel.

If $\Omega^m_{Z_\infty/R}/(\mathrm{torsion})$ is a locally free
rank $1$ module over $\mathcal{O}_{Z_\infty}$, we define the
Jacobian ideal sheaf $\mathcal{J}ac_h$ of $h$ as the annihilator
of the cokernel of the morphism
$$h^*\Omega^m_{Y_\infty/R}\rightarrow \Omega^m_{Z_\infty/R}/(\mathrm{torsion})$$
and we have $ord(Jac_h)(\psi)=ord(\mathcal{J}ac_h)(\psi)$. The
following lemma generalizes \cite[Lemma 4.1.1]{motrigid}~:

\begin{lemma}\label{jac}
Let $h:Z_\infty\rightarrow Y_\infty$ be a morphism of flat,
generically smooth $stft$ formal $R$-schemes, both equidimensional
of relative dimension $m$. Let $R'$ be a finite extension of $R$.
 For any global section $\omega$ of $\Omega^m_{Y_\eta/K}$, and any
 $\psi'\in Z_\infty(R')$,
$$ord(h^*\omega)(\psi)=ord(\omega)(h(\psi))+ord(Jac_h)(\psi)$$
\end{lemma}
\begin{proof}
This follows immediately from the definitions.
\end{proof}

\begin{lemma}\label{arcchange}
Suppose that $Y_\infty$ is a regular flat $stft$ formal
$R$-scheme, and that $Y_s$ is a strict normal crossing divisor.
 Let $\omega$ be a global section of
$\Omega^m_{Y_\eta/K}$. We denote by $\widetilde{\omega(e)}$ the
pullback of $\omega$ to the generic fiber of
$\widetilde{Y(e)_\infty}$. Let $R'$ be a finite extension of
$R(e)$, and let $\psi(e)$ be a section in
$Sm(\widetilde{Y(e)_\infty})(R')$. If we denote by $\psi$ its
image in $Y_\infty(R')$, then
$$ord(\omega)(\psi)=ord(\widetilde{\omega(e)})(\psi(e))$$
\end{lemma}
\begin{proof}
Consider the morphism $$h:Sm(\widetilde{Y(e)_\infty})\rightarrow
Y(e)_\infty$$ Since $\Omega^m_{Y(e)_\infty/R(e)}\cong
\Omega^m_{Y_\infty/R}\otimes_R R(e)$, it suffices to show,  by
Lemma \ref{jac}, that the natural map
$$(\psi(e)^*h^*\Omega^m_{\widetilde{Y(e)_\infty}/R(e)})/(\mathrm{torsion})\rightarrow
(\psi(e)^*\Omega^m_{Y(e)_\infty/R(e)})/(\mathrm{torsion})$$ is
surjective. This follows from the explicit description of the
normalization in the proof of Lemma \ref{normal}.
\end{proof}

\subsection{The order of a top form along a component of the special fiber}
We recall that $X_\infty$ denotes a regular, flat $stft$ formal
$R$-scheme, of pure relative dimension $m$, such that the special
fiber $X_s$ is a strict normal crossing divisor $\sum_{i\in I}N_i
E_i$. Let $\omega$ be a gauge form on the generic fiber $X_\eta$.
We denote by $\xi_i$ the generic point of $E_i$, for each $i$.
\begin{lemma}\label{dvr}
The local ring $\mathcal{O}_{X_\infty,\xi_i}$ of $X_\infty$ at
$\xi_i$, is a discrete valuation ring.
\end{lemma}
\begin{proof}
Locally at $\xi_i$, the divisor $E_i$ is defined by an equation
$x_0=0$. The quotient $\mathcal{O}_{X_\infty,\xi_i}/(x_0)$ is
isomorphic to $\mathcal{O}_{E_i,\xi_i}$, which is a field. Hence,
$x_0$ generates the maximal ideal of
$\mathcal{O}_{X_\infty,\xi_i}$, and since this ring is a
noetherian integral domain, it is a discrete valuation ring.
\end{proof}

\begin{lemma}\label{rank}
The $(\mathcal{O}_{X_\infty,\xi_i})$-module
$\Omega_i:=(\Omega^m_{X_\infty/R})_{\xi_i}/(\mathrm{torsion})$ is
free of rank one.
\end{lemma}
\begin{proof}
Since $(\Omega^m_{X_\infty/R})_{\xi_i}$ is finite over
$\mathcal{O}_{X_\infty,\xi_i}$, and $\mathcal{O}_{X_\infty,\xi_i}$
is a PID, the module $\Omega_i$ is free over
$\mathcal{O}_{X_\infty,\xi_i}$.

By
\cite[1.5]{formrigIII}, we have
$$\Omega^m_{X_{\infty}/R}\otimes_R K\cong \Omega^m_{X_\eta/K}$$
and since $X_\eta$ is smooth, this is a rank $1$ module over
$\mathcal{O}_{X_\eta}$. 
Hence, $\Omega_i$ has
rank $1$ over $\mathcal{O}_{X_\infty,\xi}$.
\end{proof}
%

\begin{definition}\label{order}
If $\omega\in \Omega^m_{X_\infty/R}(X_\infty)$, we define the
order of $\omega$ along $E_i$ as the length of the
$\mathcal{O}_{X_\infty,\xi_i}$-module
$\Omega_i/(\mathcal{O}_{X_\infty,\xi_i}\omega)$,
 and we denote it by $ord_{E_i}\omega$.

If $\omega\in \Omega^m_{X_\eta/K}(X_\eta)$, there exists an
integer $a\geq 0$ such that $t^a\omega\in
\Omega^m_{X_\infty/R}(X_\infty)$, by the isomorphism of sheaves
\cite[1.5]{formrigIII}
$$\Omega^m_{X_\eta/K}\cong \Omega^m_{X_\infty/R}\otimes_R K$$ and the fact
that $X_\infty$ is quasi-compact. We define the order of $\omega$
along $E_i$ as $$ord_{E_i}\omega:=ord_{E_i}(t^a\omega)-aN_i$$
\end{definition}

This definition does not depend on $a$. If $X_\infty$ is smooth,
it coincides with the one in \cite[4.3]{formner}. It follows from
Lemma \ref{dvr} and Lemma \ref{rank} that $ord_{E_i}\omega$ is
finite if $\omega$ is a gauge form on $X_\eta$.

\begin{lemma}\label{ord0}
Let $Y_\infty$ be a generically smooth $stft$ formal $R$-scheme,
endowed with an \'etale morphism of $stft$ formal $R$-schemes
$h:Y_\infty\rightarrow X_\infty$, such that the image of $h$
contains the point $\xi_i$. For any
$\omega\in\Omega^m_{X_\eta/K}(X_\eta)$, and for any connected
component $C$ of $h^{-1}(E_i)$, we have
$$ord_{E_i}\omega=ord_C h^{*}\omega.$$
\end{lemma}
\begin{proof}
Denote by $\xi'$ the generic point of $C$. Since $h$ is \'etale,
the local morphism $h^{*}:\mathcal{O}_{X_\infty,\xi_i}\rightarrow
\mathcal{O}_{Y_\infty,\xi'}$ is a flat, unramified monomorphism,
so we have isomorphisms of $\mathcal{O}_{Y_{\infty},\xi'}$-modules
\begin{eqnarray*}(\Omega^m_{X_\infty/R})_{\xi_i}
\otimes_{\mathcal{O}_{X_{\infty},\xi_i}}
\mathcal{O}_{Y_{\infty},\xi'}&\cong&
(\Omega^m_{Y_\infty/R})_{\xi'}
\\ (\Omega_i:=(\Omega^m_{X_\infty/R})_{\xi_i}/(\mathrm{torsion}))
\otimes_{\mathcal{O}_{X_{\infty},\xi_i}}
\mathcal{O}_{Y_{\infty},\xi'}&\cong&
\Omega':=(\Omega^m_{Y_\infty/R})_{\xi'}/(\mathrm{torsion})
\\
(\Omega_i/(\mathcal{O}_{X_{\infty},\xi_i}\omega))
\otimes_{\mathcal{O}_{X_{\infty},\xi_i}}
\mathcal{O}_{Y_{\infty},\xi'}&\cong&
\Omega'/(\mathcal{O}_{Y_{\infty},\xi'}h^*\omega)
\end{eqnarray*}
Now the result follows from the following algebraic property: if
$g:A\rightarrow A'$ is a flat, unramified morphism of discrete
valuation rings, if $M$ is a free $A$-module of rank $1$, and $m$
is an element of $M$, then the length of the $A$-module $M/(Am)$
equals the length of the $A'$-module $(M\otimes A')/(A'm)$.
Indeed: fixing an isomorphism of $A$-modules $A\cong M$, the
length of $A/(Am)$ is equal to the valuation of $m$ in $A$.
\end{proof}

We'll need the following technical lemma.

\begin{lemma}\label{texnic}
Let $x$ be a closed point of a generically smooth, flat, $stft$
formal $R$-scheme $Y_\infty$ of pure relative dimension $m$, and
let $\omega$ be a global section of $\Omega^m_{Y_\infty/R}$ which
induces a gauge form on $Y_\eta$. Let $\mathfrak{P}$ be a prime
ideal of the completed local ring
$\widehat{\mathcal{O}}_{Y_\infty,x}$ with $t\notin \mathfrak{P}$.
If we denote by $\widehat{\mathcal{O}}_{\mathfrak{P}}$ the
localization of $\widehat{\mathcal{O}}_{Y_\infty,x}$ at
$\mathfrak{P}$, and if we put
$$\widehat{\Omega}_{\mathfrak{P}}:=((\Omega^m_{Y_\infty/R})_x\otimes
\widehat{\mathcal{O}}_{\mathfrak{P}})/(\mathrm{torsion})$$ then
$\omega\notin \mathfrak{P}\widehat{\Omega}_{\mathfrak{P}}$.
\end{lemma}
\begin{proof}
We'll denote $\widehat{\mathcal{O}}_{Y_\infty,x}$ by
$\widehat{\mathcal{O}}$. Note that
$$(\Omega^m_{Y_\infty/R})_x\otimes\OO\cong
\Omega^m_{\mathrm{Spf}\,\OO/R}(\mathrm{Spf}\,\OO)=\Omega^m_{\OO/R}$$
Consider the tube $]x[$ of $x$ in $X_\infty$; it is an open rigid
subspace of $Y_\eta$, and coincides with the generic fiber of
$\widehat{X_\infty/x}=\mathrm{Spf}\,\OO$.

We put $C_0=\OO\otimes_R K$ and
$\Omega^m_0=\Omega^m_{\OO/R}\otimes_R K$. By
\cite[7.1.9]{dj-formal}, maximal ideals $M$ of $C_0$ are in
canonical bijective correspondence with the points $z$ of $]x[$,
and the completions of the respective local rings are isomorphic.
Moreover, by \cite[7.1.12]{dj-formal},
$$\left(\Omega^m_{]x[/K}\right)_z\otimes_{\mathcal{O}_{]x[,z}}\OO_{]x[,z}\cong \Omega^m_0\otimes_{C_0}\widehat{(C_0)_M}$$
Since $]x[$ is smooth and $\omega$ is a gauge form, this implies
that $\Omega^m_0$ is a free $C_0$-module, and $\omega$ is a
generator.

Since $\mathfrak{P}$ does not contain $t$,
$\mathfrak{Q}=\mathfrak{P}C_0$ is a prime ideal of $C_0$. Consider
the flat ring morphism $i:\OO_{\mathfrak{P}}\rightarrow
(C_0)_{\mathfrak{Q}}$ and the map of $\OO_{\mathfrak{P}}$-modules
$$(\Omega^m_{Y_\infty/R})_x\otimes
\widehat{\mathcal{O}}_{\mathfrak{P}}\rightarrow
\Omega^m_{0}\otimes (C_0)_{\mathfrak{Q}}$$ Since $\Omega^m_{0}$ is
free, and $i$ is flat, this map factors though
$$\widehat{\Omega}_{\mathfrak{P}}\rightarrow
\Omega^m_{0}\otimes (C_0)_{\mathfrak{Q}}$$ If $\omega\in
\mathfrak{P}\widehat{\Omega}_{\mathfrak{P}}$, then its image is
contained in $\mathfrak{Q}(\Omega^m_{0}\otimes
(C_0)_{\mathfrak{Q}})$, which contradicts the fact that $\omega$
is a gauge form on $Y_\eta$.
\end{proof}

For each $i\in I$, we denote by $\mathcal{I}_{E_i}$ the defining
ideal sheaf of $E_i$ in $X_\infty$. For any finite extension $R'$
of $R$ and any $\psi\in X_\infty(R')$, we denote
$ord(\mathcal{I}_{E_i})(\psi)$ by $ord_{E_i}(\psi)$ (see
Definition \ref{orderideal}). If $R'$ has ramification degree $e$
over $R$, and the closed point $\psi(0)$ of the section $\psi$ is
contained in $E_i^o$, then the equality
$t=x_i^{N_i}*\mathrm{(unit)}$ in $\mathcal{O}_{X_\infty,\psi(0)}$
implies that $ord_{E_i}(\psi)=e/N_i$.

\begin{lemma}\label{arc}
Fix a non-empty subset $J$ of $I$, and an integer $e>0$. Let $R'$
be a finite extension of $R$, of ramification index $e$, and let
$\psi$ be an element of $X_\infty(R')$, such that its closed point
$\psi(0)$ lies on
$E_J^o$. 
For any gauge form $\omega$ on $X_\eta$,
$$ord(\omega)(\psi)=\sum_{i\in J}ord_{E_i}(\psi)(ord_{E_i}\omega-1)+\max_{i\in J}\{ord_{E_i}(\psi)\}$$
\end{lemma}
\begin{proof}
 We
may assume that $\omega\in \Omega^m_{X_\infty/R}(X_\infty)$. For
notational reasons, we identify $J$ with a subset $\{0,\ldots,p\}$
of $\{0,\ldots,m\}$. We denote by $\xi_j$ the generic point of
$E_j$, for $j\in J$, and by $\xi$ the generic point of $E_J$. We
denote by $t'$ a uniformizing parameter for $R'$, with $(t')^e=t$.

By Lemma \ref{ord0}, and Lemma \ref{jac}, we can pass to an
\'etale cover and we may assume that we can find a regular system
of local parameters $(x_0,\ldots,x_m)$ on $X_\infty$ at
$x:=\psi(0)$, such that $t=\prod_{i\in J}x_i^{N_i}$.

 We can write $\omega$ as
$$\sum_{i=0}^{m}h_i(dx_0\wedge\ldots\wedge\widehat{dx_i}\wedge\ldots
\wedge dx_m),$$ with $h_i\in \mathcal{O}_{X_\infty,x}$.

Fix $j\in J$. The equality $t=\prod_{i\in J}x_{i}^{N_i}$ implies
that, for any $i\in J$,
$$dx_0\wedge\ldots\wedge\widehat{dx_i}\wedge\ldots
\wedge dx_m =(-1)^{j-i} \frac{N_i x_j}{N_j
x_i}dx_0\wedge\ldots\wedge\widehat{dx_j}\wedge\ldots \wedge dx_m$$
in $(\Omega^m_{X_\infty/R})_{\xi_j}/(\mathrm{torsion})$, so this
module is generated by
$$dx_0\wedge\ldots\wedge\widehat{dx_j}\wedge\ldots \wedge dx_m,$$
and $ord_{E_j}\omega$ equals the order of $\sum_{i\in
J}(-1)^{i}h_i\frac{N_i x_j}{x_i}$ in the discrete valuation ring
$\mathcal{O}_{X_\infty,\xi_j}$. Multiplying with the unit
$\prod_{i\in J\setminus\{j\}}x_i$, we see that $ord_{E_j}\omega$
equals the order of $$\alpha:=\sum_{i\in J}\left((-1)^{i}h_iN_i
\prod_{\ell\in J\setminus\{i\}}x_\ell\right)$$ in the discrete
valuation ring $\mathcal{O}_{X_\infty,\xi_j}$, which is the same
as the order $ord_{x_j}\alpha$ of $\alpha$ w.r.t. $x_j$ in the
unique factorization domain $\mathcal{O}_{X_\infty,x}$ (i.e. the
largest integer $N$ such that $(x_j)^N$ divides $\alpha$).

A similar computation shows that
$$ord(\omega)(\psi)=ord_{t'} \sum_{i\in J}(-1)^{i}\frac{\psi^*(h_ix_j)}{\psi^*(x_i)}N_i$$ where we chose $j\in J$ such that
$ord_{t'}\psi^*(x_j)$ is the maximum of
$\{ord_{t'}\psi^*(x_i)\,|\,i\in J\}$. Hence, it suffices to show
that
$$ord_{t'}\psi^*\alpha=\sum_{i\in J}ord_{x_i}\alpha\,.\,ord_{t'}\psi^*(x_i)$$
However,
since $\omega$ is gauge on $X_\eta$, $\alpha$ is the product of a
unit and a monomial in $(x_i)_{i\in J}$. Indeed: let
$\mathfrak{P}$ be a prime ideal in the completed local ring
$\widehat{\mathcal{O}}_{X_\infty,x}$ which does not contain $t$,
and suppose that $\alpha\in\mathfrak{P}$. Denote by
$\widehat{\mathcal{O}}_{\mathfrak{P}}$ the localization of
$\widehat{\mathcal{O}}_{X_\infty,x}$ at $\mathfrak{P}$, and put
$$\widehat{\Omega}_{\mathfrak{P}}:=((\Omega^m_{X_\infty/R})_{x}\otimes \OO_{\mathfrak{P}})/(\mathrm{torsion})$$
Since all the $x_i,i\in J$ are units in the local ring
$\widehat{\mathcal{O}}_{\mathfrak{P}}$, we see by the same
computation as above that $\alpha$ divides $\omega$ in
$\widehat{\Omega}_{\mathfrak{P}}$, and hence  $\omega \in
\mathfrak{P}\widehat{\Omega}_{\mathfrak{P}}$, which is impossible
by Lemma \ref{texnic}. So we may conclude that the only prime
divisors of $\alpha$ in the unique factorization domain
$\mathcal{O}_{X_\infty,x}$ are the prime divisors of $t$, i.e. the
elements $x_i$, modulo multiplication with a unit.

 The lemma now
follows from the fact that $ord_{t'}\psi^*x_i=ord_{E_i}(\psi)$ for
each $i\in J$.
\end{proof}

\begin{prop}\label{ord1}
Let $X_\infty$ be a regular, flat $stft$ formal $R$-scheme, such
that the special fiber $X_s$ is a strict normal crossing divisor
$\sum_{i\in I}N_i E_i$. Let $\omega$ be a gauge form on the
generic fiber $X_\eta$. Take a subset $J$ of $I$, with $|J|>1$,
and $E_J^o\neq \emptyset$. Let $h:X'_\infty\rightarrow X_\infty$
be the formal blow-up with center $E_J$, and denote by $E'_0$ its
exceptional component. We have
$$ord_{E'_0}\omega=\sum_{i\in J}ord_{E_i}\omega.$$
\end{prop}
\begin{proof}
Put $e=\sum_{i\in J}N_i$. We can find a finite unramified
extension $R'$ of $R(e)$, and a section $\psi'$ in
$X'_\infty(R')$, such that the closed point of $\psi'$ lies on
$(E'_0)^o$. Denote by $\psi$ the image of $\psi'$ in
$X_\infty(R')$. Note that $ord_{E_i}(\psi)=1$ for each $i\in J$,
since these orders are strictly positive, and $\sum_{i\in J}
ord_{E_i}(\psi)N_i$ equals $e$. Applying Lemma \ref{arc} to the
section $\psi'$, we get
$$ord(\omega)(\psi')=ord_{E'_0}\omega$$ Applying Lemma \ref{arc}
to the section $\psi$, we get
$$ord(\omega)(\psi)=\sum_{i\in J}ord_{E_i}\omega-|J|+1$$
At the generic point of $E_0'$, the module
$\Omega^m_{X'_\infty/R}/(\mathrm{torsion})$ is free of rank $1$,
and the Jacobian ideal sheaf of $h$ is the defining ideal sheaf of
$(|J|-1)E'_0$, so we may conclude by Lemma \ref{jac}.
\end{proof}

\begin{prop}\label{ord2}
Let $X_\infty$ be a regular, flat $stft$ formal $R$-scheme, of
pure relative dimension $m$, such that the special fiber $X_s$ is
a strict normal crossing divisor $\sum_{i\in I}N_i E_i$. Let
$\omega$ be a gauge form on the generic fiber $X_\eta$. Fix an
integer $e>0$. Denote by $\widetilde{\omega(e)}$ the pullback of
$\omega$ to the generic fiber of $\widetilde{X(e)_\infty}$. For
each $i\in I$, with $N_i|e$, and each connected component $C$ of
$\widetilde{E(e)}^o_i$, we have
$$ord_{C}(\widetilde{\omega(e)})= (e/N_i).ord_{E_i}\omega$$
\end{prop}
\begin{proof}
Fix $i\in I$, with $N_i|e$. Since $\widetilde{X_\infty (e)}$ is
smooth along $C$, by Lemma \ref{normal}, we can find an unramified
finite extension $R'$ of $R(e)$, and a section $\psi(e)$ in
$(\widetilde{X_\infty (e)})(R')$ such that its closed point lies
on $C$. We denote by $\psi$ the image of $\psi(e)$ in
$X_\infty(R')$. By Lemma \ref{arc},
$$ord(\omega)(\psi)=(e/N_i)ord_{E_i}\omega\mbox{ and }ord(\widetilde{\omega(e)})(\psi(e))=ord_C(\widetilde{\omega(e)})$$
Now we can apply Lemma \ref{arcchange}.
\end{proof}
\section{Computation of the volume Poincar\'e
series}\label{serreseries} Throughout this section, we put
$R=k[[t]]$, with $k$ a
 field of characteristic zero.
Let $X_\infty$ be a generically smooth $stft$ formal scheme over
$R$, with generic fiber $X_\eta$. Let $\omega$ be a gauge form on
$X_\eta$. For any integer $d>0$, we'll denote by $\omega(d)$ the
pullback of $\omega$ to $X(d)_\eta$.

\begin{definition}
For any integer $d>0$, and any locally closed subset $Z$ of $X_s$,
we put
$$F(X_\infty,\omega;d):=\int_{X(d)_\infty}|\omega(d)|\in \mathcal{M}_{X_s}$$
and we define $F_Z(X_\infty,\omega;d)$ as the image of
$F(X_\infty,\omega;d)$ under the base change morphism
$\mathcal{M}_{X_s}\rightarrow \mathcal{M}_Z$. This defines
functions
$$F(X_\infty,\omega)\,:\,\N^{*}\rightarrow \mathcal{M}_{X_s} \mbox{ and }F_Z(X_\infty,\omega)\:\,:\N^{*}\rightarrow \mathcal{M}_Z$$
which we call the local singular series (resp. local singular
series with support in $Z$) associated to the pair
$(X_\infty,\omega)$.
\end{definition}

Using terminology from \cite[4.4]{ClLo}, we define the volume
Poincar\'e series as a Mellin transform of the local singular
series.

\begin{definition}[Volume Poincar\'e series]\label{volume}
 The volume Poincar\'e series
$S(X_\infty,\omega;T)$ of the pair $(X_\infty,\omega)$ is the
generating series
$$S(X_\infty,\omega;T)=\sum_{d>0}F(X_\infty,\omega;d)T^d\in \mathcal{M}_{X_s}[[T]]$$
Its image in $\mathcal{M}_k[[T]]$ only depends on the pair
$(X_\eta,\omega)$, and is given by
$$S(X_\eta,\omega;T)=\sum_{d>0}(\int_{X(d)_\eta}|\omega(d)|)T^d$$

For any locally closed subset $Z$ of $X_s$, the volume Poincar\'e
series $S_Z(X_\infty,\omega;T)$ with support in $Z$ is defined as
the image of $S(X_\infty,\omega;T)$ under the base change morphism
$$\mathcal{M}_{X_s}[[T]]\rightarrow \mathcal{M}_{Z}[[T]]$$
\end{definition}
\begin{definition}[Serre Poincar\'e series]\label{serreser}
The Serre Poincar\'e series $S(X_\infty;T)$ of $X_\infty$ is the
generating series
$$S(X_\infty;T)=\sum_{d>0}S(X(d)_\infty)T^d\in K_0(Var_{X_s})/(\LL_{X_s}-[X_s])[[T]].$$
Its image in $K_0(Var_k)/(\LL-1)[[T]]$ only depends on $X_\eta$,
and is given by
$$S(X_\eta;T)=\sum_{d>0}S(X(d)_\eta)T^d.$$

For any locally closed subscheme $Z$ of $X_s$, the Serre
Poincar\'e series with support in $Z$ is given by
$$S_Z(X_\infty;T)=\sum_{d>0}S_Z(X(d)_\infty)T^d\in K_0(Var_{Z})/(\LL_{Z}-[Z])[[T]].$$
\end{definition}

The series $S(X_\infty,\omega;T)$ specializes to the Serre
Poincar\'e series $S(X_\infty;T)$ under the morphism
$$\mathcal{M}_{X_s}\rightarrow
\mathcal{M}_{X_s}/(\LL_{X_s}-[X_s])\cong
K_0(Var_{X_s})/(\LL_{X_s}-[X_s])$$ Likewise,
$S_Z(X_\infty,\omega;T)$ specializes to $S_Z(X_\infty;T)$

If $\mX$ is a regular, flat $stft$ formal $R$-scheme, such that
$X_s$ has strict normal crossings, and $\omega$ is a gauge form on
$X_\eta$, we can give an explicit expression for the volume
Poincar\'e series $S(X_\infty,\omega;T)$. First, we need some
technical lemmas.

\begin{lemma}\label{torus}
Let $F$ be any field, and consider a torus
$$\mathbb{G}^r_{m,F}:=\mathrm{Spec}\,F[x_1,x_1^{-1},\ldots,x_r,x_r^{-1}]$$
Take $u\in F^{\ast}$, and integers $a,b_1,\ldots,b_r$ with $a> 0$
and $gcd(a,b_1,\ldots,b_r)=1$. Then the \'etale cover
$$T:=\mathrm{Spec}\,F[x_1,x_1^{-1},\ldots,x_r,x_r^{-1}][z]/(z^a-u\prod_{j=1}^{r}x_j^{b_j})$$ is
isomorphic to $\mathbb{G}^r_{m,F}$ over $F$.
\end{lemma}
\begin{proof}
By B\'ezout, there exist integers $\alpha,\beta_1,\ldots,\beta_r$
with $\alpha.a+\sum_{j=1}^{r}\beta_j.b_j=1$. We can write $T$ as
$$T\cong \mathrm{Spec}\,F[x_1,x_1^{-1},\ldots,x_r,x_r^{-1}][z]/((u^{-\alpha}z)^a-\prod_{j=1}^{r}(u^{\beta_j}x_j)^{b_j}) $$
and hence, we may as well assume that $u=1$. In this case, writing
\begin{eqnarray*} & F[x_1,x_1^{-1},\ldots,x_r,x_r^{-1}][z]/(z^a-\prod_{j=1}^r
x_j^{b_j})\\ \cong &
F[z,z^{-1},x_1,x_1^{-1},\ldots,x_r,x_r^{-1}]/(1-z^{-a}\prod_{j=1}^r
x_j^{b_j})\end{eqnarray*} we see that $T$ is isomorphic to the
torus $\mathrm{Spec}\,F[L]$, where $L$ is the lattice
$$\Z^{r+1}/(-a.e_0+\sum_{j=1}^{r}b_j.e_j) \cong \Z^{r}$$ ($e_0,\ldots,e_r$ denotes the standard
basis of $\Z^{r+1}$). Note that $L$ is torsion-free because
$gcd(a,b_1,\ldots,b_r)=1$.
\end{proof}

\begin{lemma}\label{nochange}
Let $X_\infty$ be a regular $stft$ formal scheme over $R$, such
that its special fiber $X_s=\sum_{i=1}^{q} N_i E_i$ is a strict
normal crossing divisor. Let $J$ be a nonempty subset of
$I=\{1,\ldots,q\}$, with $|J|>1$, and let
$\pi_X:X'_\infty\rightarrow X_\infty$ be the formal blow-up with
center $E_J$. Denote its exceptional divisor by $E'_0$, and the
strict transform of $E_i$ by $E'_i$, for $i\in I$. For each subset
$K$ of $I$, the stratum $(\widetilde{E}')_{K\cup\{0\}}^{o}$ is a
Zariski piecewisely trivial fibration over $\widetilde{E}_{J\cup
K}^{o}$, and its fiber is a torus of dimension $|J\setminus K|-1$
(where we give the empty set dimension $-1$). This fibration is
compatible with the structural morphisms to $X_s$. In particular,
$[(\widetilde{E}')_{K\cup\{0\}}^{o}]=(\LL_{X_s}-[X_s])^{|J\setminus
K|-1}[\widetilde{E}_{J\cup K}^{o}]$ in $K_0(Var_{X_s})$.
\end{lemma}
\begin{proof}
Fix a subset $K$ of $I$. Covering $E_{J\cup K}^{o}\subset
X_\infty$ by opens $U$ in $X_\infty$, we may assume that
$$t=u\prod_{j\in J\cup K}x_j^{N_j},$$
with $u$ a unit, and with $x_j$ defining $E_j$. We put
$$G:=(E')^o_{K\cup\{0\}}\times_{E^o_{J\cup K}}\widetilde{E}^o_{J\cup K}$$
This means that $G$ is the \'etale cover of $(E')^o_{K\cup\{0\}}$
obtained by taking a $m_{J\cup K}$-th root of $\pi_X^*(u)^{-1}$;
we write this as $$G=(E')^o_{K\cup\{0\}}[w]/(\pi_X^*(u)w^{m_{J\cup
K}}-1)$$

We may assume that $1\in J\setminus K$. We abbreviate $J\setminus
\{1\}$ to $J^{-}$. Putting $x'_j=x_j$ for $j\in
(K\cup\{1\})\setminus J^-$, and $x'_j=x_j/x_1$, for $j\in J^-$,
 we can write
$$t=\pi_X^*(u)(x'_1)^{N_0}\!\!\!\!\prod_{j\in J^-\cup K}
(x'_j)^{N_j}$$ on $X'_\infty\setminus E'_1$, with $N_0=\sum_{i\in
J} N_i$. The locally closed subset $(E')_{K\cup\{0\}}^{o}$ is
defined by $x'_j=0$ for $j\in K\cup\{1\}$, and $x'_j\neq 0$ for
$j\in J^-\setminus K$.

We have an \'etale $(E')_{K\cup\{0\}}^{o}$-morphism
\begin{eqnarray*}(\widetilde{E}')^{o}_{K\cup
\{0\}}=(E')_{K\cup\{0\}}^{o}[z]/(\pi_X^*(u)\!\prod_{j\in
J^-\setminus K}(x'_j)^{N_j}z^{m_{K\cup \{0\}}}-1)
\\ \longrightarrow G=
(E')_{K\cup\{0\}}^{o}[w]/(\pi_X^*(u)w^{m_{J\cup K}}-1),
\end{eqnarray*}
 defined by
 $$w\mapsto z^{m_{K\cup \{0\}}/m_{J\cup K}}\prod_{j\in
J^-\setminus K}(x'_j)^{N_j/m_{J\cup K}}$$ We see it is induced by
the isomorphism
$$(\widetilde{E}')^{o}_{K\cup
\{0\}}=G[z]/(z^{{m_{K\cup \{0\}}/m_{J\cup K}}}-w\prod_{j\in
J^-\setminus K}(x'_j)^{-N_j/m_{J\cup K}})$$

We claim that the fiber $E_{x}$ of $(\widetilde{E}')^{o}_{K\cup
\{0\}}$ over any point $x$ of $\widetilde{E}^o_{J\cup K}$ is
isomorphic to $\G_{m,k(x)}^{|J\setminus K|-1}$. This concludes the
proof, by \cite{sebag1}, Theorem 4.2.3.

So let us prove our claim. Denote by $G_x$ the fiber of $G$ over
$x$. The regular function $w$ is constant along $G_x$, with value
$w(x)\in k(x)^{\times}$. By \cite[IV-26]{Eisenbud}, the
restrictions of the functions $x'_j$, $j\in J^-\setminus K$, to
$G_x$, induce an isomorphism $G_x\cong
\mathrm{Spec}\,k(x)[\,x'_j,(x'_j)^{-1}]_{j\in J^-\setminus K}$ (if
$A$ is an $R$-algebra topologically of finite type, the formal
blow-up of $\mathrm{Spf}\,A$ at a regular center $(x_j)_{j\in J}$
contained in the special fiber, is the formal completion of the
blow-up of $\mathrm{Spec}\,A$ at the center $(x_j)_{j\in J}$; in
particular, the special fibers are the same; so
\cite[IV-26]{Eisenbud} carries over to the formal case).

 The fiber
$E_x$ is isomorphic to
$$G_x[z]/(z^{m_{K\cup \{0\}}/m_{J\cup
K}}-w(x)\prod_{j\in J^-\setminus K}(x'_j)^{-N_j/m_{J\cup K}})$$ We
can conclude by Lemma \ref{torus} and the fact that
$$gcd(m_{K\cup\{0\}},N_j)_{j\in J^-\setminus K}= m_{J\cup K}$$
\end{proof}

\begin{theorem}\label{expression}
Let $X_\infty$ be a regular flat $stft$ formal scheme over $R$, of
pure relative dimension $m$, such that $X_s=\sum_{i\in I} N_i E_i$
is a strict normal crossing divisor. Let $\omega$ be a gauge form
on $X_\eta$.

 If we denote by $\mu_i$ the order
$ord_{E_i}\omega$ of $\omega$ along $E_i$, then, for any integer
$d>0$,
$$F(X_\infty,\omega;d)=\LL^{-m}\sum_{\emptyset\neq J\subset
I}(\LL-1)^{|J|-1}[\widetilde{E}_J^{o}](\sum_{\stackrel{k_i\geq
1,i\in J}{\sum_{i\in J} k_iN_i=d}}\LL^{-\sum_i k_i \mu_i}\,)\
\mbox{in}\ \mathcal{M}_{X_s}$$ In particular,
$$S(X(d)_\infty)=\sum_{i\in I, N_i|d}[\widetilde{E}_i^{o}]\ \mbox{in}\
K_0(Var_{X_s})/(\LL_{X_s}-[X_s])$$
\end{theorem}
\begin{proof}
First, suppose that $d$ is not $X_s$-linear. By Theorem
\ref{neronmodel}, $Sm(\widetilde{X(d)_\infty})\rightarrow
X(d)_\infty$ is a N\'eron smoothening. By Lemma \ref{ord2}, for
each $i\in I$ with $N_i|d$, the order of $\widetilde{\omega(d)}$
along any component of $\widetilde{E(d)}_i^o$ equals $d\mu_i/N_i$.

By \cite{motrigid}, Prop. 4.3.1 and \cite{NiSe2}, Lemma 6.4, we
see that
\begin{eqnarray*}
F(X_\infty,\omega;d)&=&\LL^{-m}\sum_{N_i|d}[\widetilde{E}_i^{o}]\LL^{-
d\mu_i/N_i}
\\&=&\LL^{-m}\sum_{\emptyset\neq J\subset
I}(\LL-1)^{|J|-1}[\widetilde{E}_J^{o}](\sum_{\stackrel{k_i\geq
1,i\in J}{\sum_{i\in J} k_iN_i=d}}\LL^{-\sum_i
k_i\mu_i}\,)\,,\qquad\qquad(*)\end{eqnarray*}
 provided that
$d$ is not $X_s$-linear.

We use Lemma \ref{nochange} to show that the expression (*) does
not change under the formal blow-up $h$ of a stratum $E_J$. Denote
the exceptional divisor of $h$ by $E'_0$, and the strict transform
of $E_i$ by $E'_i$, for $i\in I$.

Observe that the multiplicity $N_0$ of $E'_0$ in $X'_s$ equals
$\sum_{j\in J}N_j$, and that, by Lemma \ref{ord1},
$\mu_0:=ord_{E'_0}\omega=\sum_{j\in J}\mu_j$. For any non-empty
subset $L$ of $I$, containing $J$, and any vector $k=(k_i)$ in
$(\N^{*})^{L}$, denote by $M(k)$ the set of indices $j\in J$ such
that $k_j=m(k):=\min\{k_i\,|\,i\in J\}$.
 We get, for each subset $L$ of $I$
containing $J$,
\begin{eqnarray*}
&&(\LL-1)^{|L|-1}[\widetilde{E}_L^{o}](\sum_{\stackrel{k_i\geq
1,i\in L}{\sum_{i\in L} k_iN_i=d}}\LL^{-\sum_i k_i \mu_i}\,) \\&=&
\sum_{L\setminus J\subset K\subset
L}(\LL-1)^{|L|-1}[\widetilde{E}_{L}^{o}](\sum_{\stackrel{k=(k_i)\in
(\N^{*})^{L}, M(k)=L\setminus K}{\sum_{i\in L}
k_iN_i=d}}\LL^{-\sum_i k_i\mu_i}\,)
\\&=&\sum_{L\setminus J\subset K\subset L}(\LL-1)^{|K|}[(\widetilde{E}')_{K\cup\{0\}}^{o}]
\\&& (\sum_{\stackrel{k=(k_i)\in (\N^{*})^{L}}{M(k)=L\setminus
K,\,\gamma(k)}}\LL^{-m(k)(\sum_{i\in J}\mu_i)-\sum_{i\in J}
(k_i-m(k))\mu_i-\sum_{i\in L\setminus J}k_i\mu_i}\,)
\end{eqnarray*}
where condition $\gamma(k)$ means $$m(k)(\sum_{i\in J}
N_i)+\sum_{i\in J} (k_i-m(k))N_i+\sum_{i\in L\setminus
J}k_iN_i=d.$$ We can rewrite this last expression as
$$\sum_{K\subset L, J\cup
K=L}(\LL-1)^{|K|}[(\widetilde{E}')_{K\cup\{0\}}^{o}](\sum_{\stackrel{k'=(k'_i)\in
(\N^{*})^{K\cup\{0\}}}{\sum_{i}k'_iN_i=d}}\LL^{-\sum_ik'_i\mu_i}\,)\,.$$

By Lemma \ref{linear}, we can find, for each $d>0$, a composition
of blow-ups $\pi^{(j)}:X_\infty^{(j+1)}\rightarrow
X_\infty^{(j)}$, $j=0,\ldots,r-1$, such that
\begin{itemize}
\item $X_\infty^{(0)}=X_\infty$, \item the special fiber of
$X^{(j)}_\infty$ is a strict normal crossing divisor
$$X_s^{(j)}=\sum_{i=1}^{k_j} N_i^{(j)}E_i^{(j)},$$ \item $\pi^{(j)}$
is a blow-up with center $E_{J^{(j)}}$, for some non-empty subset
$J^{(j)}$ of $\{1,\ldots,k_j\}$, \item $d$ is not
$X_s^{(r)}$-linear.
\end{itemize}
Since the expression $(*)$ does not change under any of the
$\pi^{(j)}$, it is valid for all $d$.
\end{proof}

\begin{corollary}\label{serrat}
Let $X_\infty$ be a generically smooth $stft$ formal $R$-scheme,
of pure relative dimension $m$. Suppose that $X_\infty$ admits a
resolution $X'_\infty\rightarrow X_\infty$, with special fiber
$X'_s=\sum_{i\in I} N_i E_i$. Let $\omega$ be a gauge form on
$X_\eta$.

 The volume
Poincar\'e series $S(X_\infty,\omega;T)$ is rational over
$\mathcal{M}_{X_s}$. In fact, if we put $\mu_i:=ord_{E_i}\omega$,
then the series is given explicitly by
$$S(X_\infty,\omega;T)=\LL^{-m}\sum_{\emptyset\neq J\subset I}(\LL-1)^{|J|-1}[\widetilde{E}_J^{o}]
\prod_{i\in J}\frac{\LL^{-\mu_i}T^{N_i}}{1-\LL^{-\mu_i}T^{N_i}} \
\mathrm{in}\ \mathcal{M}_{X_s}[[T]]$$ In particular,
$$S(X_\infty;T)=\sum_{i\in I}[\widetilde{E}^o_i]\frac{T^{N_i}}{1-T^{N_i}}\ \in\ \left(K_0(Var_{X_s})/(\LL_{X_s}-[X_s])\right)[[T]]$$
\end{corollary}
\begin{proof}
We may assume that $X'_\infty=X_\infty$.  By Theorem
\ref{expression},
\begin{eqnarray*}
S(X_\infty,\omega;T)&=&\LL^{-m}\sum_{d>0}\sum_{\emptyset\neq
J\subset
I}\left((\LL-1)^{|J|-1}[\widetilde{E}_J^{o}](\sum_{\stackrel{k_i\geq
1,i\in J}{\sum_{i\in J} k_iN_i=d}}\LL^{-\sum_i
k_i\mu_i}\,)\right)T^d
\\&=&\LL^{-m}\sum_{\emptyset\neq J\subset
I}\left((\LL-1)^{|J|-1}[\widetilde{E}_J^{o}]\prod_{i\in
J}\sum_{k_i>0}\LL^{-k_i\mu_i}T^{k_iN_i}\right)
\\&=&\LL^{-m}\sum_{\emptyset\neq J\subset
I}\left((\LL-1)^{|J|-1}[\widetilde{E}_J^{o}]\prod_{i\in
J}\frac{\LL^{-\mu_i}T^{N_i}}{1-\LL^{-\mu_i}T^{N_i}}\right)
\end{eqnarray*}
\end{proof}


By Proposition \ref{affineres}, any affine generically smooth flat
$stft$ formal $R$-scheme admits a resolution of singularities. By
the additivity of the motivic integral, we obtain an expression
for the volume Poincar\'e series in terms of a finite atlas of
local resolutions. In particular, we obtain the following result.

\begin{cor}\label{rationality}
Let $X_\infty$ be a generically smooth $stft$ formal $R$-scheme,
of pure relative dimension $m$. Let $\omega$ be a gauge form on
$X_\eta$.
 The volume
Poincar\'e series $S(X_\infty,\omega;T)$ is rational over
$\mathcal{M}_{X_s}$. More precisely, there exists a finite subset
$S$ of $\Z\times \N^{\ast}$ such that $S(X_\infty,\omega;T)$
belongs to the subring
$$\mathcal{M}_{X_s}\left[\frac{\LL^{a}T^b}{1-\LL^{a}T^b} \right]_{(a,b)\in S} $$
of $\mathcal{M}_{X_s}[[T]]$.
\end{cor}

\section{The motivic volume}\label{tame}
Suppose that $R=k[[t]]$, with $k$ a field of characteristic zero.

Let $X_\infty$ be a generically smooth, flat $stft$ formal
$R$-scheme, of pure relative dimension $m$, and let $\omega$ be a
gauge form on $X_\eta$.
 It is not
possible to associate a motivic Serre invariant to
$X_\infty\widehat{\times}\widehat{R^s}$ in a direct way, since the
normalization $R^s$ of $R$ in $K^s$ is not a discrete valuation
ring. We will define a motivic object by taking a limit of Serre
invariants of finite ramifications of $X_\infty$, instead.

\begin{definition}[\cite{GLM}, (2.8)]
There is a unique $\mathcal{M}_{X_s}$-linear morphism
$$\lim_{T\to\infty}:\mathcal{M}_{X_s}\left[\frac{\LL^{a}T^b}{1-\LL^{a}T^b} \right]_{(a,b)\in \Z\times \N^{\ast}}\longrightarrow \mathcal{M}_{X_s}$$
mapping $$\prod_{(a,b)\in I}\frac{\LL^{a}T^b}{1-\LL^{a}T^b}$$ to
$(-1)^{|I|}=(-1)^{|I|}[X_s]$, for each finite subset $I$ of
$\Z\times \N^{\ast}$. We call the image of an element its limit
for $T\to\infty$.
\end{definition}
To see that this morphism is well-defined, note that the image of
an element is given by the constant term of its Taylor development
in $T^{-1}$.
\begin{prop}\label{separable}
 The
limit of $-S(X_\infty,\omega;T)$ for $T\to \infty$ is
well-defined, and does not depend on the choice of $\omega$. If
$\mX'\rightarrow \mX$ is any resolution of singularities, with
$X'_s=\sum_{i\in I}N_i E_i$, then this limit is given explicitly
by
$$\LL^{-m}\sum_{\emptyset\neq J\subset
I}(1-\LL)^{|J|-1}[\widetilde{E}_J^{o}]$$ in $\mathcal{M}_{X_s}$.
\end{prop}
\begin{proof}
This follows immediately from the computation in Corollary
\ref{serrat}, and the remark preceding Corollary
\ref{rationality}.
\end{proof}

\begin{definition}\label{motvolume}
 The
motivic volume
$$S(X_\infty;\widehat{K^s})\in
\mathcal{M}_{X_s}$$ is by definition the limit of
$-S(X_\infty,\omega;T)$ for $T\to\infty$, where $\omega$ is any
gauge form on $X_\eta$.

When $V$ is a locally closed subset of $X_s$, we define the
motivic volume $S_{V}(X_\infty;\widehat{K^s})$ with support in $V$
as the image of $S(X_\infty;\widehat{K^s})$ under the base change
morphism
$$\mathcal{M}_{X_s}\rightarrow \mathcal{M}_{V}$$

Finally, we define the motivic volume
$$S(X_\eta;\widehat{K^s})\in\mathcal{M}_k$$ as
the image of $S(X_\infty;\widehat{K^s})$ under the forgetful
morphism $\mathcal{M}_{X_s}\rightarrow \mathcal{M}_k$. It only
depends on $X_\eta$, and not on $X_\infty$.
\end{definition}

\begin{remark}\label{global}
For these definitions, it is not necessary that $X_\eta$ admits a
global gauge form. By additivity of the motivic integral w.r.t.
open covers of $X_\infty$, we can use a system of local gauge
forms on $X_\eta$, whose domains cover $X_\eta$. Since the limit
of $-S(X_\infty,\omega;T)$ for $T\to\infty$ does not depend on
$\omega$, these local gauge forms do not have to coincide on the
intersections of their domains. See also \cite{NiSe2}, Remark 6.5.
\end{remark}

Ayoub constructs in \cite{ayoub-lecture} a triangulated category
of motives for rigid varieties $X_\eta$ over $k((t))$, and he
shows that this category is equivalent to a certain subcategory of
the stable homotopy-category of schemes over the torus
$\mathbb{G}_{m,k}$. Pull-back via the unit section of
$\mathbb{G}_{m,k}$ yields a motive over $k$, and this construction
can be used to describe Ayoub's motivic nearby cycle functor. It
seems very plausible that Ayoub's motive for $X_\eta$ coincides
with our motivic volume $S(X_\eta;\widehat{K^s})$ in an
appropriate Grothendieck ring of motives.
 Moreover, any reasonably
defined additive invariant of $\overline{X}_\eta$ (e.g. its Hodge
polynomial) should be computable on $S(X_\eta;\widehat{K^s})$.

 A priori, $S(X_\infty;\widehat{K^s})$ depends not
only on $\mX\widehat{\times}_R\widehat{R^s}$, but also on
$X_\infty$. The following Proposition shows that, for any proper
subscheme $Z$ of $X_s$,
$\chi_{top}(S_{Z}(X_\infty;\widehat{K^s}))$ depends only on the
$\widehat{K^s}$-analytic space $\overline{]Z[}$, and not on the
embedding of this space in $\mX\widehat{\times}_R\widehat{R^s}$,
if we assume that $\mX$ is algebrizable.

\begin{prop}\label{eulintr}
Suppose $k$ is algebraically closed. Let $X$ be a generically
smooth flat $R$-variety. Let $Z$ be a proper subvariety of the
special fiber $X_s$. Then
$$\chi_{top}(S_{Z}(\widehat{X};\widehat{K^s}))=\chi_{\acute{e}t}(\overline{]Z[})\,,$$
where $\chi_{\acute{e}t}$ is the Euler characteristic associated
to Berkovich' \'etale $\ell$-adic cohomology for non-archimedean
analytic spaces.
\end{prop}
\begin{proof}
Let $\varphi$ be a topological generator of the absolute Galois
group $G(K^s/K)$. By definition,
$$\chi_{top}(S_{Z}(\widehat{X};\widehat{K^s}))=-\lim_{T\to\infty}\chi_{top}(S_{Z}(\widehat{X};T))$$
Hence, by our Trace Formula in Theorem \ref{trace},
$$\chi_{top}(S_{Z}(\widehat{X};\widehat{K^s}))=-\lim_{T\to\infty}\sum_{n>0}Tr(\varphi^n\,|\,H(\overline{]Z[})\,)T^n$$

Recall the identity \cite[1.5.3]{weil1}
\begin{eqnarray*}\label{weil1}
\sum_{n>0}Tr(F^n,V)T^n&=&T\frac{d}{dT}log(det(1-TF,V)^{-1})
\\&=&-\frac{T\frac{d}{dT}(det(1-TF,V))}{det(1-TF,V)}
\end{eqnarray*} for any endomorphism $F$ on a finite dimensional vector
space $V$. Taking limits yields
$$
-\lim_{T\to\infty}\sum_{n>0}Tr(F^n,V)T^n=dim(V).$$ Applying this
to $F=\varphi$ and $V=H(\overline{]Z[})$ yields
$$\chi_{top}(S_{Z}(\widehat{X},\widehat{K^s}))=\chi_{\acute{e}t}(\overline{]Z[})\,.$$
\end{proof}

\section{Applications to motivic zeta functions and the Milnor
fibration}\label{app}
\subsection{The analytic Milnor
fiber}\label{milnor} Let $X_\infty$ be a flat $stft$ formal
$R$-scheme. For any closed point $x$ on the special fiber $X_s$,
the tube $]x[$ of $x$ in $X_\infty$ is an open rigid subspace of
$X_{\eta}$. We will call it the analytic Milnor fiber of
$X_\infty$ at $x$, and denote it by $\cF_x$. As a $K$-rigid space,
it is canonically isomorphic to the generic fiber of the
completion $\widehat{X_\infty/x}$ of $X_\infty$ along $x$ (see
\cite[0.2.7]{bert}). If $R=k[[t]]$, and $X_\infty$ is the $t$-adic
completion a morphism of $k$-varieties $f:Z\rightarrow
\mathrm{Spec}\,k[t]$, we will also call $\cF_x$ the analytic
Milnor fiber of $f$ at $x$.

\subsubsection{Cohomology of the analytic Milnor
fiber}\label{cohomology} The following proposition shows that
$\cF_x$ has the ``right" \'etale cohomology.

\begin{lemma}\label{cohmil}
Suppose that $k$ is algebraically closed. Let $X$ be a flat
$R$-variety, and let $\mathcal{F}_x$ be the analytic Milnor fiber
of $\widehat{X}$ at a closed point $x$ of $X_s$. Let $F$ be a
constructible \'etale torsion sheaf on the generic fiber
$X\times_R K$ of the $R$-scheme $X$, with torsion orders prime to
$p$. Let $\hat{F}$ be the induced sheaf on the analytic space
$X_{\eta}\widehat{\times}_K\widehat{K^{t}}$. The \'etale
cohomology space
$H^q(\cF_x\widehat{\times}\widehat{K^{t}},\hat{F})$ is canonically
isomorphic to the $q$-th cohomology space $R^q\psi_{\eta}^t(F)_x$
of the stalk
 at $x$ of the complex of tame nearby
cycles $R\psi_{\eta}^t(F)$ of $X/R$, and this isomorphism is
compatible with the geometric monodromy action of $G(K^t/K)$.
\end{lemma}
\begin{proof}
 This follows from the comparison result in
\cite[3.5]{berk-vanish2}.
%
%
\end{proof}
\begin{theorem}\label{milnorfiber}
Let $X$ be a smooth, irreducible $\C$-variety, and consider a
dominant morphism $f:X\rightarrow \A^1_{\C}=\mathrm{Spec}\,\C[t]$.
Let $x$ be a point of the hypersurface $X_s$ defined by $f$,
denote by $F_x$ the canonical topological Milnor fiber of $f$ at
$x$, and by $M$ the monodromy automorphism on the singular
cohomology $H_{sing}(F_x,\C)$ \footnote{We recalled the
definitions in the introduction.}.

Let $\cF_x$ be the analytic Milnor fiber of $f$ at $x$, and let
$\varphi$ be the canonical topological generator of
$G(\C((t))^s/\C((t))\,)\cong \hat{\Z}(1)(\C)$. Fix an embedding of
$\Q_{\ell}$ in $\C$. There are canonical isomorphisms
$$H^i_{sing}(F_x,\C)\cong H^i(\cF_x\widehat{\times}\widehat{\C((t))^s},\Q_\ell)\otimes_{\Q_{\ell}}\C$$
compatible with the action of $M$ and $\varphi$.
\end{theorem}
\begin{proof}
By Lemma \ref{cohmil}, we have, for each integer $n>0$, canonical
$G(\C((t))^s/\C((t))\,)$-equivariant isomorphisms
$$H^i(\cF_x\widehat{\times}\widehat{\C((t))^s},\Z/\ell^n)\cong R^i\psi^t_{\eta}(\Z/\ell^n)_x$$
By \cite[XIV]{SGA7b}, there exist canonical isomorphisms
$$R^i\psi^t_{\eta}(\Z/\ell ^n)_x\cong
R^i\psi^{an}_{\eta}(\Z/\ell^n)_x$$ where $R\psi^{an}_{\eta}$ is
the complex of analytic vanishing cycles, and these isomorphisms
are compatible with the action of
$G(\C((t))^s/\C((t))\,)=\hat{\Z}(1)(\C)$ and
$\pi_1(\G_{m,\C})=\Z$. Finally, there are canonical isomorphisms
$$R^i\psi^{an}_{\eta}(\Z/\ell^n)_x\cong H^i_{sing}(F_x,\Z/\ell^n)$$
respecting the action of $\pi_1(\G_{m,\C})$ (see e.g.
\cite{Kuli},(8.11.7)). Taking projective limits over $n$ yields
\begin{eqnarray*}
H^i(\cF_x\widehat{\times}\widehat{\C((t))^s},\Z_{\ell})&:=&\underleftarrow{\lim}H^i(\cF_x\widehat{\times}\widehat{\C((t))^s},\Z/\ell^n)
\\&\cong&
\underleftarrow{\lim}H^i_{sing}(F_x,\Z/\ell^n)\\&=&H^i_{sing}(F_x,\Z_{\ell})\end{eqnarray*}
Tensoring with $\C$ yields the required result.
\end{proof}


As a corollary, we recover some classical results concernig the
cohomology of the Milnor fiber. Let $X$ be a smooth irreducible
variety over $\C$, and let $f:X\rightarrow
\A^1_{\C}=\mathrm{Spec}\,\C[t]$ be a dominant morphism. Let $x$ be
a complex point of the hypersurface $X_s$ defined by $f$. Let
$h:X'\rightarrow X$ be an embedded resolution of singularities for
$f$, with $(f\circ h)=\sum_{i\in I} N_i E_i$. Let $F_x$ be the
topological Milnor fiber at $x$, let $M$ be the monodromy
transformation on the singular cohomology $H_{sing}(F_x,\C)$, and
denote by $\zeta_x(T)$ the monodromy zeta function at $x$ (i.e.
the alternating product of the characteristic polynomials of $M$).

\begin{corollary}[\cite{Acampo2}]
If $x$ is a singular point of $f$, then
$Tr(M\,|\,H_{sing}(F_x,\C))=0$. Else,
$Tr(M\,|\,H_{sing}(F_x,\C))=1$.
\end{corollary}
\begin{proof}
If $Y$ is any regular, flat $R$-variety, then
$Sm(\widehat{Y})\rightarrow \widehat{Y}$ is a N\'eron
$R$-smoothening, by \cite[\S 3, Prop 2]{neron}. Hence,
$S_x(\widehat{X})=0$ if $x$ is a singular point of $f$, and
$S_x(\widehat{X})=1$ else. Now our trace formula in Theorem
\ref{trace} yields the result.
\end{proof}
\begin{corollary}[A'Campo's formula \cite{A'C}]
\begin{eqnarray*}
Tr(M^d\,|\,H_{sing}(F_x,\C))&=&\sum_{N_i|d}N_i\chi_{top}(E_i^o\cap
h^{-1}(x)),
\\ \zeta_x(T)&=&\prod_i(T^{N_i}-1)^{-\chi_{top}(E^o_i\cap h^{-1}(x))}.
\end{eqnarray*}
\end{corollary}
\begin{proof}
The second formula follows from the first via the identity
(\ref{weil1}) in Section \ref{tame}.

By Theorem \ref{trace}, the expressions are valid if we replace
$F_x$ by $\cF_x\widehat{\times}\widehat{\C((t))^s},$ where $\cF_x$
is the analytic Milnor fiber $\cF_x$ at $x$, and if we replace
singular cohomology by Berkovich' \'etale cohomology for analytic
spaces.

The expression now follows from the comparison result in Theorem
\ref{milnorfiber}.
\end{proof}

\subsubsection{Points of the analytic Milnor fiber}\label{points}
We suppose that $R=k[[t]]$, with $k$ an algebraically closed field
of characteristic zero.  Let $X$ be any variety over $k$, and
consider a $k$-morphism $f:X\rightarrow \mathrm{Spec}\,k[t]$, flat
over the origin. Denote by $\widehat{X}$ the $t$-adic completion
of $f$. We can describe the points of the generic fiber
$X_{\eta}$, and the specialization map $sp:|X_{\eta}|\rightarrow
|\widehat{X}|=|X_s|$ on the level of the underlying sets.

We denote by $\mathcal{L}(X)$ the arc scheme of $X$, as defined in
\cite[p.1]{DLinvent}. It is a scheme over $k$, of infinite type in
general, and for any field $k'$ containing $k$,
$$\mathcal{L}(X)(k')=Hom_k(\Spec k'[[u]],X)$$ Reduction modulo $u$ yields a
projection morphism of $k$-schemes
$\pi_0:\mathcal{L}(X)\rightarrow X$.

For any integer $d>0$, we denote by $\mathcal{X}(d)$ the closed
subscheme of $\mathcal{L}(X)$ defined by
$$\mathcal{X}(d)=\{\psi\in\mathcal{L}(X)\,|\,f(\psi)=u^d\}$$
We will construct a canonical bijection
$$\varphi:X_{\eta}(K(d))\rightarrow \mathcal{X}(d)(k)$$
such that the square

$$\begin{CD}
X_{\eta}(K(d))@>\varphi>> \mathcal{X}(d)(k)
\\ @VspVV @VV\pi_0V
\\ X_s(k)@>=>> X_s(k)
\end{CD}$$
commutes.

The specialization morphism of ringed sites $sp:X_\eta\rightarrow
\widehat{X}$ induces a bijection $X_\eta(K(d))\rightarrow
\widehat{X}(R(d))$, and the morphism $sp:X_\eta(K(d))\rightarrow
X_s(k)$ maps a point of $X_\eta(K(d))$ to the reduction modulo
$t_d$ of the corresponding point of $\widehat{X}(R(d))$.

 By Grothendieck's
Existence Theorem \cite[5.4.1]{ega3}, the completion functor
induces a bijection $(X\times_{k[t]} R)(R(d))\rightarrow
\widehat{X}(R(d))$. A reparametrization $t\mapsto t_d$ yields a
bijection $(X\times_{k[t]} R)(R(d))\rightarrow \mathcal{X}(d)(k)$.

Hence, for any closed point $x$ of $X_s$, the set $\cF_x(K(d))$
corresponds bijectively to the fiber of
$\pi_0:\mathcal{X}(d)(k)\rightarrow X_s(k)$ over $x$. In other
words, a $K(d)$-point of $\cF_x$ is nothing but an arc $\psi\in
Hom_k(\mathrm{Spec}\,k[[u]], X)$ with $f(\psi(u))=u^d$ and
$\psi(0)=x$.

  The Galois group
$G(K(d)/K)=\mu_d(k)$ acts on $X_\eta(K(d))$ as follows: if $\psi$
is an element of $\mathcal{X}(d)(k)$, an element $\xi\in\mu_d(k)$
acts on $\psi$ by $\xi*\psi(u)=\psi(\xi\cdot u)$.

\subsection{The local singular series and the Gelfand-Leray
form}\label{singseries} Let $k$ be a field of characteristic zero,
let $X$ be a smooth irreducible variety over $k$, of dimension
$m$, and let $f:X\rightarrow \A^{1}_k=\mathrm{Spec}\,k[t]$ be a
dominant morphism, smooth on the complement of the special fiber
$X_s$ of $f$ in $X$.
 Denote by
$\widehat{X}$ the $t$-adic completion of $X$, and by $X_\eta$ its
generic fiber. Let $\omega$ be a gauge form on $X$. The Koszul
complex associated to $df$ is exact on $V:=X\setminus X_s$, and in
particular
$$\begin{CD}
\Omega^{m-2}_V@>\wedge df>>\Omega^{m-1}_V@>\wedge
df>>\Omega^{m}_V@>>> 0
\end{CD}$$
is an exact sequence of sheaves on $V$. This means that we can
choose, for each point $y$ in $V$, an element $\alpha_y$ in
$(\Omega^{m-1}_V)_y$, such that $(df)_y\wedge \alpha_y=\omega_y$.
Since the $\alpha_y$ are unique up to a term $df\wedge \!*\,$,
they glue together to a well-defined global section
$$\frac{\omega}{df}\in \Omega^{m-1}_{V/\A^1_k}(V)$$
This relative differential form induces a gauge form
$\frac{\omega}{df}$ on $X_\eta$.

\begin{definition}\label{localsing0}
 We call $$\frac{\omega}{df}\in \Omega^{m-1}_{X_\eta/K}(X_\eta)$$ the Gelfand-Leray
form of $\omega$ w.r.t. $f$.
\end{definition}
%

\begin{lemma}\label{order-gl}
Let $h:X'\rightarrow X$ be an embedded resolution of singularities
for the morphism $f$, with $X'_s=\sum_{i\in I}N_i E_i$ and
relative canonical divisor $K_{X'/X}=\sum_{i\in I}(\nu_i-1)E_i$.
Then
$$ord_{E_i}h^*(\frac{\omega}{df})=\nu_i - N_i$$ for any $i\in I$.
\end{lemma}
\begin{proof}
The pullback $h^{*}(\omega/df)$ coincides with the Gelfand-Leray
form associated to $h^{*}(\omega)$ and $f\circ h$, since for any
differential form $\alpha$ on $X$, we have
$$h^{*}\alpha\wedge d(f\circ h)=h^{*}\alpha\wedge h^{*}df=h^{*}(\alpha\wedge df).$$

 Choosing local
coordinates $(x_1,\ldots,x_m)$ and $(x'_1,\ldots,x'_m)$ on $X$,
resp. $X'$, we may assume that $f\circ h=u(x'_1)^{N_i}$, with $u$
a unit, and that the Jacobian of $h$ is given by
$Jac_h=u'(x_1)^{\nu_i-1}$, for some unit $u'$. This means that the
pullback $h^{*}\omega$ is given by $(x_1)^{\nu_i-1}$ times a gauge
form on $X'$. Passing to an \'etale cover of $X'$, we may assume
that $u=1$, by Lemma \ref{ord0}. In this case,
$h^{*}\omega/d(f\circ h)$ equals $(x'_1)^{\nu_i-N_i}$ times a
nowhere vanishing relative form on $\widehat{X'}/R$, hence
$ord_{E_i}(h^{*}\omega/d(f\circ h))=\nu_i- N_i$.
\end{proof}

\begin{cor}\label{localseries} We have, for any $d>0$,
$$F(\widehat{X},\frac{\omega}{df};d)=\LL^{-(m-1)}\sum_{\emptyset\neq J\subset I}
(\LL-1)^{|J|-1}[\widetilde{E}^o_J]\,\left( \sum_{\stackrel{k_i\geq
1,i\in J}{\sum_{i\in J} k_iN_i=d}}\LL^{k_i(N_i-\nu_i)}\right)\,\in
\mathcal{M}_{X_s}$$ In particular, the right hand side does not
depend on the choice of $\omega$.

If $d$ is not $X'_s$-linear, then the above expression reduces to
$$F(\widehat{X},\frac{\omega}{df};d)=\LL^{-(m-1)}\sum_{N_i|d}
[\widetilde{E}^o_i]\,\LL^{d-d\nu_i/N_i}\,\in \mathcal{M}_{X_s}$$
\end{cor}
\begin{proof}
This follows immediately from Theorem \ref{expression}.
\end{proof}

\begin{definition}\label{localsing}
We define the local singular series associated to $f$ by
$$F(f;d):=F(\widehat{X},\frac{\omega}{df};d)$$
where $\omega$ is any gauge form on $X$.
\end{definition}
\begin{remark}
For this definition, it is not necessary that $X$ admits a global
gauge form $\omega$. See Remark \ref{global}.
\end{remark}
\subsection{Comparison to the motivic zeta
function}\label{comparzeta} Let $k$ be a field of characteristic
zero, let $X$ be a smooth irreducible $k$-variety of dimension
$m$, and let $f:X\rightarrow \A^1_k=\mathrm{Spec}\,k[t]$ be a
dominant morphism.

 As in \cite[p.1]{DLinvent}, we denote, for any integer
$d>0$, by $\mathcal{L}_d(X)$ the $k$-scheme representing the
functor
$$(k-algebras)\rightarrow (Sets)\,:\,A\mapsto Hom_k(\Spec A[u]/(u^{d+1}),X)$$
Following \cite[3.2]{DL3}, we denote by $\mathcal{X}_{d,1}$ the
$X_s$-variety
$$\mathcal{X}_{d,1}=\{\psi\in
\mathcal{L}_{d}(X)\,|\,f(\psi(u))=u^d\,mod\,u^{d+1}\}$$ where the
structural morphism $\mathcal{X}_{d,1}\rightarrow X_s$ is given by
reduction modulo $u$. In \cite[3.2.1]{DL3}, the motivic zeta
function $Z(T)$ of $f$ is defined as
$$Z(T)=\sum_{d=1}^{\infty}[\mathcal{X}_{d,1}]\LL^{-md}T^{d}\in \mathcal{M}_{X_s}[[T]]$$
Actually, the coefficients live in a more refined Grothendieck
ring, as we will see in Section \ref{motmilnor}.

 In \textit{loc. cit.}, Denef and Loeser
show that $Z(T)$ is rational over $\mathcal{M}_{X_s}$, and they
give an explicit expression in terms of an embedded resolution of
$f$ (see also Section \ref{motmilnor}).

\begin{lemma}\label{weilzeta}
Let $X$ be a smooth irreducible variety over $k$ of dimension $m$,
 let $f:X\rightarrow \A^1_k$ be a
 dominant morphism, and let $\omega$ be a gauge form on $X$. We
 have
$$F(\widehat{X},\frac{\omega}{df};d)=\LL^{-(d+1)(m-1)}[\mathcal{X}_{d,1}]$$ in $\mathcal{M}_{X_s}$, for each integer
$d>0$.
\end{lemma}
\begin{proof}
Take an embedded resolution $X'\rightarrow X$ for $f$, such that
$d$ is not $X'_s$-linear, and combine Corollary \ref{localseries}
with the computation in \cite{DLLefschetz}, Theorem 2.4.
\end{proof}

\begin{theorem}\label{weil}
Let $X$ be a smooth irreducible variety over $k$ of dimension $m$,
 let $f:X\rightarrow \A^1_k$ be a
 dominant morphism, and let $\omega$ be any gauge form on $X$. We
 have
$$S(\widehat{X},\frac{\omega}{df};T)=\LL^{-(m-1)}Z(\LL T)\in
\mathcal{M}_{X_s}[[T]]$$
\end{theorem}
\begin{proof}
This follows immediately from Lemma \ref{weilzeta}.
\end{proof}

Hence, we recover the motivic zeta function as a Mellin transform
of the local singular series associated to $f$. Our trace formula
yields the following cohomological interpretation of the motivic
zeta function.

\begin{theorem}\label{dl} Suppose that $k$ is algebraically
closed.
 Let $Z$ be a subvariety of $X_s$, proper
over $k$, and denote by $\mathcal{X}_{d,1,Z}$ the fiber of
$\mathcal{X}_{d,1}\rightarrow X_s$ over $Z$. Let $\varphi$ be a
topological generator of the geometric monodromy group $G(K^s/K)$.
 For any integer $d>0$,
$$\chi_{top}(\mathcal{X}_{d,1,Z})=Tr(\varphi^{d}|H(\overline{]Z[}))\,.$$
\end{theorem}
\begin{proof}
This is an immediate consequence of Lemma \ref{weilzeta}, and the
trace formula Theorem \ref{trace}.
\end{proof}
We recover the main result of Denef and Loeser's paper
\cite{DLLefschetz}.
\begin{cor}[\cite{DLLefschetz}, Theorem 1.1]
If $k=\C$, we get, for each complex point $x$ on $X_s$, and each
$d>0$,
$$\chi_{top}(\mathcal{X}_{d,1,x})=Tr(M^d\,|\,H_{sing}(F_x,\C))\,,$$
where $F_x$ is the canonical topological Milnor fiber of $f$ at
$x$, and $M$ is the monodromy transformation.
\end{cor}
\begin{proof}
This follows from the comparison result in Theorem
\ref{milnorfiber}, and from Theorem \ref{dl} above.
\end{proof}

\subsection{The motivic Milnor fiber}\label{motmilnor}
We recall Denef and Loeser's definition of the motivic Milnor
fiber \cite[3.5]{DL3}. Let $k$ be an algebraically closed field of
characteristic zero, let $X$ be a smooth irreducible $k$-variety
of dimension $m$, and let $f:X\rightarrow
\A^1_k=\mathrm{Spec}\,k[t]$ be a dominant morphism. We denote the
hypersurface defined by $f$ in $X$ by $X_s$.

Let $\mu$ be the inverse limit of the groups $\mu_d(k)$ of $d$-th
roots of unity in $k$. Let $S$ be a variety over $k$. By a good
$\mu$-action on an algebraic variety $Z$ over $S$, we mean an
action of $\mu$ on $Z$,
 equivariant w.r.t. the
structure morphism $Z\rightarrow S$ (where $S$ carries the trivial
action), that factors through some $\mu_d$, and such that each
orbit
 is contained in an affine open subscheme of $Z$. The
relative equivariant Grothendieck ring $K_0^{\mu}(Var_S)$ of
$S$-varieties with good $\mu$-action is defined in
\cite[2.4]{DL3}.

For each $d>0$, the $X_s$-variety $\mathcal{X}_{d,1}$ carries a
good $\mu$-action, defined as follows: if $\xi$ is an element of
$\mu_d(k)$, and $\psi(u)$ is an element of $\mathcal{X}_{d,1}$,
then $\xi*\psi=\psi(\xi\cdot u)$. Denef and Loeser define the
equivariant motivic zeta function as
$$Z^{mon}(T)=\sum_{d=1}^{\infty}[\mathcal{X}_{d,1}]\LL^{-md}T^d \in K_0^{\mu}(Var_{X_s}).$$
They show it is rational over
$\mathcal{M}^{\mu}_{X_s}=K_0^{\mu}(Var_{X_s})[\LL_{X_s}^{-1}]$,
where $\LL_{X_s}$ denotes the class of the affine line over $X_s$,
with trivial $\mu$-action.

In fact, they obtain an explicit formula in terms of an embedded
resolution of singularities for $f$ on $X$. Let $h:X'\rightarrow
X$ be an embedded resolution for $f$, with $X'_s=\sum_{i\in I}N_i
E_i$, and with relative canonical divisor $K_{X'|X}=\sum_{i\in
I}(\nu_i-1)E_i$. For any non-empty subset $J$ of $I$, the cover
$\widetilde{E}_J^o/E_J^o$ (defined in Section \ref{weakner}) is
Galois with group $\mu_{m_J}(k)$.

By \cite{DL3}, Theorem 3.3.1, we have
$$Z^{mon}(T)=\sum_{\emptyset\neq J\subset I}(\LL-1)^{|J|-1}[\widetilde{E}_J^{o}]
\prod_{i\in J}\frac{\LL^{-\nu_i}T^{N_i}}{1-\LL^{-\nu_i}T^{N_i}} \
\mathrm{in}\ \mathcal{M}^{\mu}_{X_s}[[T]].$$ Inspired by the
$p$-adic case \cite{Denef5}, Denef and Loeser defined the motivic
 nearby cycles $\mathcal{S}_f$ by taking formally the
limit of $-Z^{mon}(T)$ for $T\to\infty$, i.e.
$$\mathcal{S}_f= \sum_{\emptyset\neq J\subset I}(1-\LL)^{|J|-1}[\widetilde{E}_J^{o}]\in \mathcal{M}^{\mu}_{X_s}.$$
This terminology is justified by the fact that, when $k=\C$, for
each closed point $x$ of $X_s$, the mixed Hodge structure of the
fiber $\mathcal{S}_{f,x}\in \mathcal{M}^{\mu}_{k}$ of
$\mathcal{S}_f$ over $x$ coincides with the mixed Hodge structure
of the Milnor fiber of $f$ at $x$ (in an appropriate Grothendieck
group of mixed Hodge structures), and the $\mu$-action corresponds
to the semi-simple part of the monodromy action \cite[4.2]{DL5}.
Denef and Loeser called $\mathcal{S}_{f,x}$ the motivic Milnor
fiber of $f$ at $x$.

\begin{theorem}
Denoting the $t$-adic completion of $X$ by $\widehat{X}$, we have
$$S(\widehat{X};\widehat{K^s})=\LL^{-(m-1)}\mathcal{S}_f \in\mathcal{M}_{X_s}$$
For any closed point $x$ on $X_s$, we have
$$S_x(\widehat{X};\widehat{K^s})=\LL^{-(m-1)}\mathcal{S}_{f,x} \in\mathcal{M}_{x}$$
\end{theorem}
\begin{proof}
This follows from Theorem \ref{weil}.
\end{proof}

Hence, we recover Denef and Loeser's motivic nearby cycles as the
motivic volume of the rigid ``nearby fiber'' $X_\eta$, and in some
sense, we recover the motivic Milnor fiber $\mathcal{S}_{f,x}$ as
the motivic volume of the analytic Milnor fiber $\mathcal{F}_x$.

\section{Acknowledgements} The authors would like to thank
Fran\c{c}ois Loeser for inspiring discussions. The first author is
indebted to L. Illusie, V. Berkovich, and O. Villamayor, who were
so kind to answer his questions on nearby cycles, non-archimedean
spaces, and resolution of singularities.
\bibliographystyle{hplain}
\bibliography{wanbib,wanbib2}

\begin{thebibliography}{10}

\bibitem{sga7a}
{\em Groupes de monodromie en g\'eom\'etrie alg\'ebrique. {I}}.
\newblock Springer-Verlag, Berlin, 1972.
\newblock S\'eminaire de G\'eom\'etrie Alg\'ebrique du Bois-Marie 1967--1969
  (SGA 7 {I}), Dirig\'e par A. Grothendieck. Avec la collaboration de M.
  Raynaud et D. S. Rim, Lecture Notes in Mathematics, Vol. 288.

\bibitem{SGA7b}
{\em Groupes de monodromie en g\'eom\'etrie alg\'ebrique. {II}}.
\newblock Springer-Verlag, Berlin, 1973.
\newblock S\'eminaire de G\'eom\'etrie Alg\'ebrique du Bois-Marie 1967--1969
  (SGA 7 II), Dirig\'e par P. Deligne et N. Katz, Lecture Notes in Mathematics,
  Vol. 340.

\bibitem{Acampo2}
N.~A'Campo.
\newblock Le nombre de {L}efschetz d'une monodromie.
\newblock {\em Indag. Math.}, 35:113--118, 1973.

\bibitem{A'C}
N.~A'Campo.
\newblock La fonction z\^eta d'une monodromie.
\newblock {\em Comment. Math. Helvetici}, 50:233--248, 1975.

\bibitem{ayoub-lecture}
J.~Ayoub.
\newblock Motives of rigid varieties.
\newblock {\em Lecture at the Conference on Hodge Theory, Venice, 19--24 June
  2006}.

\bibitem{Berk-etale}
V.~G. Berkovich.
\newblock {\'Etale cohomology for non-Archimedean analytic spaces.}
\newblock {\em Publ. Math., Inst. Hautes \'Etud. Sci.}, 78:5--171, 1993.

\bibitem{berk-vanish2}
V.~G. Berkovich.
\newblock {Vanishing cycles for formal schemes II.}
\newblock {\em Invent. Math.}, 125(2):367--390, 1996.

\bibitem{bert}
P.~Berthelot.
\newblock {Cohomologie rigide et cohomologie rigide \`{a} supports propres}.
\newblock {\em Prepublication, Inst. Math. de Rennes}, 1996.

\bibitem{bosch}
S.~Bosch.
\newblock {\em {Lectures on formal and rigid geometry.}}
\newblock preprint, http://www.math1.uni-muenster.de/sfb/about/publ/bosch.html,
  2005.

\bibitem{BGR}
S.~Bosch, U.~G{\"{u}}ntzer, and R.~Remmert.
\newblock {\em {Non-Archimedean analysis. A systematic approach to rigid
  analytic geometry.}}, volume 261 of {\em {Grundlehren der Mathematischen
  Wissenschaften}}.
\newblock Springer Verlag, 1984.

\bibitem{neron}
S.~Bosch, W.~{L\"u}tkebohmert, and M.~Raynaud.
\newblock {\em {N\'eron models}}, volume~21.
\newblock {Ergebnisse der Mathematik und ihrer Grenzgebiete}, 1990.

\bibitem{formrigIII}
S.~Bosch, W.~{L\"u}tkebohmert, and M.~Raynaud.
\newblock {Formal and rigid geometry. III: The relative maximum principle.}
\newblock {\em Math. Ann.}, 302(1):1--29, 1995.

\bibitem{formner}
S.~Bosch and K.~Schl{\"o}ter.
\newblock N\'eron models in the setting of formal and rigid geometry.
\newblock {\em Math. Ann.}, 301(2):339--362, 1995.

\bibitem{VMayor}
A.~Bravo, S.~Encinas, and O.~Villamayor.
\newblock A simplified proof of desingularization and applications.
\newblock {\em Rev. Mat. Iberoamericana}, 21(2):349–--458, 2005.

\bibitem{ClLo}
R.~Cluckers and F.~Loeser.
\newblock Constructible motivic functions, integrals with parameters, and cell
  decomposition.
\newblock arxiv:math.AG/0410203.

\bibitem{conrad}
B.~Conrad.
\newblock Irreducible components of rigid spaces.
\newblock {\em Ann. Inst. Fourier}, 49(2):473--541, 1999.

\bibitem{dj-formal}
A.~J. de~Jong.
\newblock {Crystalline {D}ieudonn\'e module theory via formal and rigid
  geometry.}
\newblock {\em Publ. Math., Inst. Hautes \'Etud. Sci.}, 82:5--96, 1995.

\bibitem{weil1}
P.~Deligne.
\newblock La conjecture de {W}eil. {I}.
\newblock {\em Publ. Math., Inst. Hautes \'Etud. Sci.}, 43:273--307, 1973.

\bibitem{DenefBour}
J.~Denef.
\newblock Report on {I}gusa's local zeta function.
\newblock In {\em S\'eminaire Bourbaki, Vol. 1990/91, Exp. No.730-744}, volume
  201-203, pages 359--386, 1991.

\bibitem{Denef5}
J.~Denef.
\newblock Degree of local zeta functions and monodromy.
\newblock {\em Compositio Math.}, 89:207--216, 1993.

\bibitem{DL5}
J.~Denef and F.~Loeser.
\newblock Motivic {I}gusa zeta functions.
\newblock {\em J. Algebraic Geom.}, 7:505--537, 1998, arxiv:math.AG/9803040.

\bibitem{DLinvent}
J.~Denef and F.~Loeser.
\newblock Germs of arcs on singular algebraic varieties and motivic
  integration.
\newblock {\em Invent. Math.}, 135:201--232, 1999, arxiv:math.AG/9803039.

\bibitem{DL3}
J.~Denef and F.~Loeser.
\newblock Geometry on arc spaces of algebraic varieties.
\newblock {\em Progr. Math.}, 201:327--348, 2001, arxiv:math.AG/0006050.

\bibitem{DLLefschetz}
J.~Denef and F.~Loeser.
\newblock {Lefschetz numbers of iterates of the monodromy and truncated arcs.}
\newblock {\em Topology}, 41(5):1031--1040, 2002.

\bibitem{Eisenbud}
D.~Eisenbud and J.~Harris.
\newblock {\em The geometry of schemes}, volume 197 of {\em Graduate Texts in
  Mathematics}.
\newblock Springer-Verlag, New York, 2000.

\bibitem{ega3}
A.~Grothendieck and J.~Dieudonn\'e.
\newblock El\'ements de {G}\'eom\'etrie {A}lg\'ebrique, iii.
\newblock {\em Publ. Math., Inst. Hautes \'Etud. Sci.}, 11:5--167, 1961.

\bibitem{ega4.1}
A.~Grothendieck and J.~Dieudonn\'e.
\newblock El\'ements de {G}\'eom\'etrie {A}lg\'ebrique, iv, premi\`ere partie.
\newblock {\em Publ. Math., Inst. Hautes \'Etud. Sci.}, 20:5--259, 1964.

\bibitem{GLM}
G.~Guibert, F.~Loeser, and M.~Merle.
\newblock {Iterated vanishing cycles, convolution, and a motivic analogue of a
  conjecture of {S}teenbrink.}
\newblock {\em Duke Math. J.}, 132(3):409--457, 2006.

\bibitem{Igusa:intro}
J.~Igusa.
\newblock {\em An introduction to the theory of local zeta functions}.
\newblock Studies in {A}dvanced {M}athematics. AMS, 2000.

\bibitem{Illusie}
L.~Illusie.
\newblock Th\'eorie de {B}rauer et caract\'eristique d'{E}uler-{P}oincar\'e
  (d'apr\`es {P}. {D}eligne).
\newblock In {\em The Euler-Poincar\'e characteristic}, volume~82 of {\em
  Ast\'erisque}, pages 161--172. 1981.

\bibitem{Kempf}
G.~Kempf, F.~Knudsen, D.~Mumford, and B.~Saint-Donat.
\newblock {\em Toroidal embeddings 1}, volume 339 of {\em Lecture Notes in
  Mathematics}.
\newblock Springer-Verlag, 1973.

\bibitem{Ko}
M.~Kontsevich.
\newblock Lecture at {O}rsay.
\newblock December 7, 1995.

\bibitem{Kuli}
V.~Kulikov.
\newblock {\em Mixed Hodge Structures and Singularities}, volume 132 of {\em
  Cambridge Tracts in Mathematics}.
\newblock Cambridge University Press, 1998.

\bibitem{motrigid}
F.~Loeser and J.~Sebag.
\newblock {Motivic integration on smooth rigid varieties and invariants of
  degenerations}.
\newblock {\em Duke Math. J.}, 119:315--344, 2003.

\bibitem{Milnor}
J.~Milnor.
\newblock {\em Singular points of complex hypersurfaces}, volume~61 of {\em
  Annals of Math. Studies}.
\newblock Princeton University Press, 1968.

\bibitem{NiSe3}
J.~Nicaise and J.~Sebag.
\newblock Rigid geometry and the monodromy conjecture.
\newblock {\em to appear in the Conference Proceedings of
  \textit{Singularit\'es}, CIRM, Luminy, 24 January-25 February 2005}.

\bibitem{NiSe-Milnor}
J.~Nicaise and J.~Sebag.
\newblock Invariant de {S}erre et fibre de {M}ilnor analytique.
\newblock {\em C.R.Ac.Sci.}, 341(1):21--24, 2005.

\bibitem{NiSe-weilres}
J.~Nicaise and J.~Sebag.
\newblock Motivic {S}erre invariants and {W}eil restriction.
\newblock {\em preprint}, 2006.

\bibitem{NiSe2}
J.~Nicaise and J.~Sebag.
\newblock Motivic {S}erre invariants of curves.
\newblock {\em to appear in Manuscr. Math.}, 2006.

\bibitem{Poo}
B.~Poonen.
\newblock The {G}rothendieck ring of varieties is not a domain.
\newblock {\em Math. Res. Letters}, 9(4):493--498, 2002, arxiv:math.AG/0204306.

\bibitem{RoVe}
B.~Rodrigues and W.~Veys.
\newblock Poles of {Z}eta functions on normal surfaces.
\newblock {\em Proc. London Math. Soc.}, 87(3):164--196, 2003.

\bibitem{sebag1}
J.~Sebag.
\newblock Int\'egration motivique sur les sch\'emas formels.
\newblock {\em Bull. Soc. Math. France}, 132(1):1--54, 2004.

\end{thebibliography}
\end{document}